\documentclass{amsart}

\usepackage{arxiv_math}

\title{Immersed Curves and 4-Manifold Invariants}

\author[J.Cohen]{Jesse Cohen}
\address{Fachbereich Mathematik (AZ)\\Universit\"{a}t Hamburg\\
	Bundesstra\ss{}e 55\\
	20146 Hamburg, Germany}
\email{jesse.cohen@uni-hamburg.de}
\author[G.Guth]{Gary Guth}
\address{Department of Mathematics\\Stanford University\\
	Building 380\\
	Stanford, California 94305}
\email{gmguth@stanford.edu}

\begin{document}
	\begin{abstract}
		For 3-manifolds with torus boundary, the bordered Heegaard Floer invariants of Lipshitz--Ozsv\'ath--Thurston have a geometric interpretation as immersed multi-curves with local systems in the punctured torus according to the work of Hanselman--Rasmussen--Watson. We consider morphisms between these immersed curve invariants and show that they compute certain cobordism maps. More precisely, we relate composition in the Fukaya category of immersed curves in the punctured torus to composition of morphisms between the bordered Floer invariants, which have interpretations in terms of certain cobordism maps. We make use of this formalism to obstruct smooth equivalences between 4-manifolds with boundary, and between surfaces with boundary in the 4-ball.
	\end{abstract}
	\maketitle

	\section{Introduction}
Heegaard Floer homology is a suite of invariants of 3-manifolds introduced by Ozsv\'ath and Szab\'o \cite{os_holodisks, os_holotri,os_knotinvts}. In its simplest form, Heegaard Floer homology associates to a (based, oriented, closed) three-manifold $Y$ a vector space $\HFh(Y)$ which arises as the homology of a chain complex $\CFh(Y)$ whose homotopy type is a diffeomorphism invariant of $Y$. Moreover, Heegaard Floer homology is functorial with respect to cobordisms: if $W$ is a cobordism from $Y_1$ to $Y_2$, there is an induced map $\widehat{F}_W:\HFh(Y_1)\ra\HFh(Y_2)$, which is itself a smooth invariant of the cobordism up to homotopy.

Bordered Floer homology is an extension of Heegaard Floer homology to 3-manifolds with parametrized boundary, developed by Lipshitz--Ozsv\'ath--Thurston \cite{LOT_bordered_HF}. Given a connected surface $F$, they define a differential graded algebra $\cA(F)$ and, to a 3-manifold $M$ with boundary identified with $F$, they associate a pair of invariants: a left differential graded module $\CFDh(M)$ and right $A_\infty$-module $\CFAh(M)$, whose homotopy types are smooth invariants of $M$. If $M$ and $N$ are bordered 3-manifolds and $h:\partial N \xra{\cong}\partial M$ is an orientation \emph{reversing} diffeomorphism, there is a pairing theorem \cite[Theorem 1.3]{LOT_bordered_HF} which states that 
\begin{align*}
	\CFAh(M) \boxtimes \CFDh(N) \simeq \CFh(M \cup_h N),
\end{align*}
where $\CFAh(M) \boxtimes \CFDh(N)$ is a convenient model for the derived tensor product. There is also a $\Hom$ pairing theorem \cite[Theorem 1]{LOT_HF_as_morphism}, which states that if $h: \partial N  \ra \partial M$ is an orientation \emph{preserving} diffeomorphism, then there is a homotopy equivalence
\begin{align*}
	\Mor^{\mathcal{A}(F)}(\CFDh(M),\CFDh(N)) \simeq \CFh(-M \cup_h N),
\end{align*}
where $\Mor^{\mathcal{A}(F)}(\CFDh(M),\CFDh(N))$ is the chain complex of left $\mathcal{A}(F)$-module homomorphisms $$f:\CFDh(M)\to\CFDh(N)$$ with differential given by $\partial f:=\partial_N\circ f+f\circ\partial_M$. There is an analogous statement for morphism complexes for the right $A_\infty$-module invariants. Since we will work exclusively in the case that $F=T^2$, we will suppress the boundary surface from our notation and declare $\mathcal{A}\defeq\mathcal{A}(T^2)$ throughout.

In the case of 3-manifolds with torus boundary, bordered Floer homology has a striking geometric interpretation due to Hanselman--Rasmussen--Watson \cite{hanselman2023bordered}. They show that the type D structure associated to a 3-manifold $M$ with torus boundary can viewed as an immersed multi-curve with local systems in $\partial M \smallsetminus \pt.$ The pairing theorem takes a particularly nice form as well in this theory: if $M$ and $N$ are 3-manifolds with torus boundary with immersed curve invariants $\btheta_M$ and $\btheta_N$, and $h: \partial N \ra \partial M$ is an orientation reversing diffeomorphism, then the Heegaard Floer homology of $M\cup_h N$ is homotopy equivalent to the Lagrangian intersection Floer homology
\begin{align*}
	\CFh(M \cup_h N) \simeq\hom_{\Fuk_{\circ}}(\btheta_M,h(\btheta_N))
\end{align*}
of $\btheta_M$ and $h(\btheta_N)$ in the punctured torus. 

These two results can be reformulated as a statement about the partially wrapped Fukaya category $\Fuk_{\circ}$ of the punctured torus with stop consisting of a single point $z$ in its boundary: the immersed curve invariant corresponding to $M$ is a compactly supported object of this Fukaya category, and the Heegaard Floer homology of $M \cup_h N$ is given by the hom-pairing in $\Fuk_{\circ}$, i.e. $\CFh(M \cup N)$ has a model given by $\F \langle \btheta_M \cap \btheta_N \rangle$ with differential given by counting holomorphic bigons between the two curves. 

\begin{theorem}[Hanselman--Rasmussen--Watson, {\cite{hanselman2023bordered}}]
	Suppose that $M$ and $N$ are 3-manifolds with parametrized torus boundary. Then, the bordered invariants of $M$ and $N$ can be represented by immersed curves $\btheta_M$ and $\btheta_N$ in $\partial M$ and $\partial N$ respectively. Moreover, if $h: \partial N \ra \partial M$ is an orientation preserving diffeomorphism, then there is a homotopy equivalence of chain complexes: 
    \begin{align*}
        \CFh(-M\cup_h N)\simeq \Mor_\cA(\CFAh(M), \CFAh(N))\simeq \hom_{\Fuk_{\circ}}(\btheta_M, h(\btheta_N)).
    \end{align*}
\end{theorem}

\begin{rem}
    We note that in \cite{hanselman2023bordered}, the pairing theorem takes a slightly different form: $M$ and $N$ are glued by an orientation \emph{reversing} diffeomorphism, and $\btheta_M$ is paired with $\overline{h}(\btheta_N)$, where $h$ denotes the composition of the gluing map with the elliptic involution. This pairing computes the box tensor product $\CFAh(M)\boxtimes \CFDh(N) \simeq \CFh(M \cup_{\overline{h}} N).$ In the morphism formulation, the elliptic involution plays no role. For more details, see \Cref{sec:background}. 
\end{rem}

Composition in the Fukaya category is given by counting holomorphic triangles between three sets of curves. In light of the morphism pairing theorem for the bordered invariants, one would naturally expect this identification of $\Mor_\cA(\CFAh(M_i),\CFAh(M_j))$ with $\hom_{\Fuk_{\circ}}(\btheta_1, \btheta_2)$ to be compatible with composition. Our first result demonstrates that this is indeed the case.

\begin{theorem}\label{thm:FukayaHomPairing}
	Let $M_i$ be 3-manifolds with torus boundary and let $\phi_{i+1}: \partial M_{i+1} \ra \partial M_i$ be orientation preserving diffeomorphisms. Let $\btheta_i$ be their corresponding immersed curve invariants. Fix a parametrization of $\partial M_k$; the gluing maps then determine parametrizations for the remaining $\partial M_i$. Then, there are homotopy equivalences
	\begin{align}
		\Mor_\cA(M_i,M_{i+1})\simeq\hom_{\Fuk_{\circ}}(\btheta_i,\phi_{i+1}(\btheta_{i+1})),
	\end{align}
	where $\Mor_\cA(M_i,M_j):=\Mor_\cA(\CFAh(M_i),\CFAh(M_j))$ such that the diagrams
	\begin{align}
		\begin{tikzcd}[ampersand replacement=\&]
			\Mor_\cA(M_0,M_1)\otimes\cdots\otimes\Mor_\cA(M_{k-1},M_k)\arrow[d,"\simeq"']\arrow[r,"\circ_k"] \& \Mor_\cA(M_0,M_k)\arrow[d,"\simeq"]\\
			\hom_{\Fuk_{\circ}}(\btheta_0,\phi_{ 1}(\btheta_1))\otimes\cdots\otimes\hom_{\Fuk_{\circ}}(\phi_{1, \dots, k-1}(\btheta_{k-1}),\phi_{1, \dots, k}(\btheta_k))\arrow[r,"\circ_k"] \& \hom_{\Fuk_{\circ}}(\btheta_0,\phi_{1, \dots, k}(\btheta_k))
		\end{tikzcd}
	\end{align}
	all homotopy commute. Here, $\phi_{1, \dots, j}$ is shorthand for the composition $\phi_1 \circ \cdots \circ\phi_j$.
\end{theorem}

When the gluing map is understood, we will often drop it from the notation. 

In the context of Heegaard Floer theory, maps between Floer complexes ought to be induced by 4-dimensional cobordisms. Indeed, by work of the first author \cite{CohenComposition}, the composition map
\begin{align*}
	\Mor_\cA(M_0,M_1)\otimes\Mor_\cA(M_{1},M_2)\xra{\circ_2} \Mor_\cA(M_0,M_2)
\end{align*}
is realized by a particular cobordism 
\begin{align*}
	\cP:(-M_0 \cup M_1) \amalg (-M_1\cup M_2) \ra (-M_0\cup M_2) 
\end{align*}
often referred to as the \emph{pair of pants} cobordism, which is central to the construction of cobordism maps in Heegaard Floer theory. Combined with \Cref{thm:FukayaHomPairing}, this implies that counting triangles is equivalent to computing certain cobordism maps on Heegaard Floer homology.

\begin{example}\label{ex:946_disks}
	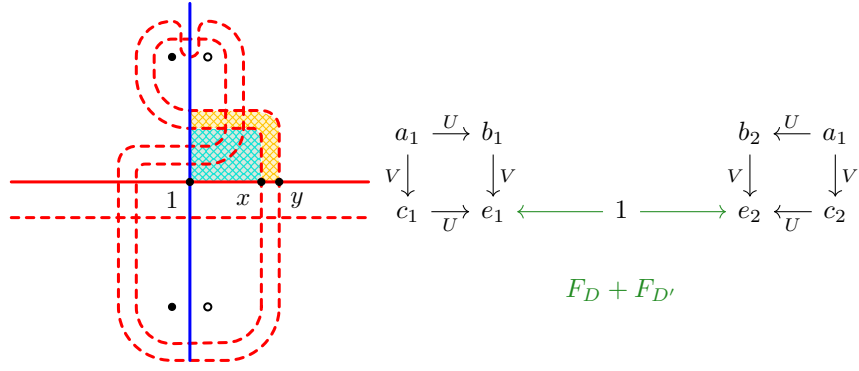
\begin{figure}[h]
		\centering
		\raisebox{-2cm}{\includegraphics[scale=0.95]{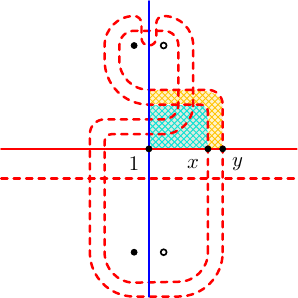}}
		\begin{tikzcd}[column sep=0.5cm,row sep=0.5cm]
			a_1\ar[r,"U"]\ar[d,"V"'] & b_1\ar[d,"V"] &  & b_2\ar[d,"V"'] & a_1\ar[l,"U"']\ar[d,"V"]\\
			c_1\ar[r,"U"'] & e_1 & 1\ar[l,ForestGreen]\ar[r,ForestGreen] & e_2 & c_2\ar[l,"U"]\\ &  & {\color{ForestGreen}F_D+F_{D'}} & &
		\end{tikzcd}
		\caption{The immersed multicurve invariant for $K = m9_{46}$ and the associated knot Floer complex. The map $F_D + F_{D'}$ is given by the small blue triangles.} 
		\label{fig:9_46_disks}
	\end{figure}
	To bring matters down to earth, consider the following motivating example. Let  $K = m9_{46}$. By the work of Hanselman-Rasmussen-Watson, the bordered module $\CFDh(\KC)$ can be realized as a compactly supported element of the Fukaya category, i.e. as an immersed multi-curve $\btheta_K$ in the torus. A lift of $\btheta_K$ to the infinite cylinder, $S^1 \times \R$, is shown in \Cref{fig:9_46_disks}. To compute the knot Floer complex, $\CFK(S^3, K)$, we may algebraically pair the invariants $\CFDh(\KC)$ and $\CFDh(S^1 \times D^2, S^1 \times \pt)$ (i.e., compute the space of morphism from the former to the latter). We denote by $\bbeta_{\id}$ the curve invariant for the solid torus containing the knot given by an $S^1$-fiber is the vertical curve in \Cref{fig:9_46_disks}. We may realize this pairing geometrically as the complex generated by the intersection points of these two curves. The differential can then be computed by counting bigons weighted by a $U$ or $V$ power given by their algebraic intersection with the two basepoints. The resulting complex is shown on the right of \Cref{fig:9_46_disks}. 
	
	The knot $K$ bounds a pair of slice disks $D$ and $D'$ which are shown in Figure \ref{fig:9_46_disk4}.
	\begin{figure}
		\centering
		\includegraphics[width=0.5\linewidth]{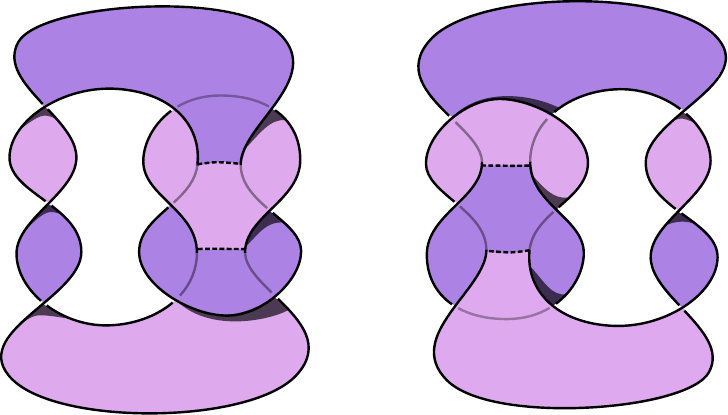}
		\caption{A pair of slice disks for $K=m9_{46}$.}
		\label{fig:9_46_disk4}
	\end{figure}
	The maps associated to these disks can be distinguished by comparing their induced maps, 
	\begin{align*}
		F_D+F_{D'}: \CFK(U) \ra \CFK(K).
	\end{align*}
	In the basis chosen above, $(F_D+F_{D'})(1) = e_1 + e_2.$ There is an associated map 
	\begin{align*}
		F_{B^4 \smallsetminus D} + F_{B^4 \smallsetminus D'}: \CFDh(S^3 \smallsetminus U) \ra \CFDh(S^3 \smallsetminus K), 
	\end{align*}
	which is determined up to homotopy by the property that $\bI\boxtimes (F_{B^4 \smallsetminus D} + F_{B^4 \smallsetminus D'}) \simeq F_D+F_{D'}$. As the notation indicates, $F_{B^4 \smallsetminus D} + F_{B^4 \smallsetminus D'}$ should be thought of as a map associated to the complements of the disks $D$ and $D'$, though we do not claim that such maps are invariants in general. Under the identification 
	\begin{align*}
		\Mor^\cA(\CFDh(S^3 \smallsetminus U), \CFDh(S^3 \smallsetminus K)) \simeq \hom_{\Fuk_{\circ}}(\btheta_U, \btheta_K),
	\end{align*}
	the map $(F_{B^4 \smallsetminus D} + F_{B^4 \smallsetminus D'})$ is represented by some collection of intersection points $\bm b$ between $\btheta_U$ and $\btheta_K$.
	
	To compute the map $(F_D + F_{D'})$ in terms of the immersed curve invariants, we are to compute the composition in the Fukaya category of the morphisms (intersection points) corresponding to 1 and $\bm b$, i.e. we count triangles with one vertex on $1 \in \CFK(U) \simeq \F\langle \btheta_U \cap \bbeta_{\id} \rangle$, the second vertex on $\bm b \in \F\langle \btheta_U \cap \btheta_K \rangle$, and the third on $e_1 + e_2 \in \CFK(K) \simeq \F\langle \btheta_K \cap \bbeta_{\id} \rangle$. Since we \emph{know} that this triangle counting map should compute $(F_D+F_{D'})(1) = e_1 + e_2,$ we may deduce that $\bm b$ must be given by the two intersection points $\bm b$ shown in \Cref{fig:9_46_disks}; there are only two intersection points which complete triangles from $1$ to $e_1 + e_2$ and have the correct boundary conditions. 
	
	While it may be philosophically gratifying to realize the cobordism map $(F_{B^4 \smallsetminus D} + F_{B^4 \smallsetminus D'})$ as some linear combination of intersection points between the curve invariants for $U$ and $K$, its utility seems rather limited, since we made use of the fact that we could already compute the map $F_D+F_{D'}$. While this approach may, at first, appear somewhat backwards, the advantage is that now the map $(F_{B^4 \smallsetminus D} + F_{B^4 \smallsetminus D'})$ has been pinned down, this map can be used to compute any other examples built by gluing some 4-manifold with corners to the complements of these disks. For instance, consider the 4-manifolds $B^4_{1/n}(D)$ and $B^4_{1/n}(D')$ obtained by performing surgery on these disks (i.e. the cobordisms obtained by capping off the surgered concordances $S_{1/n}(U) \ra S_{1/n}(K)$). These cobordism maps can easily be calculated: we simply pair $(\btheta_U\cup \btheta_K, \bm b)$ not with $\bbeta_{\id}$, but instead with the curve invariant for the $1/n$-framed solid torus, which is given by the line of slope $1/n$, and counting triangles between these three curves with a corner at $\bm b$. In particular, it is easy to see that for any $n > 0$ these two cobordisms are distinguished by their induced maps. See \Cref{fig:disk_surgery} for an example.
	
	\begin{figure}
		\centering
		\includegraphics[scale=1]{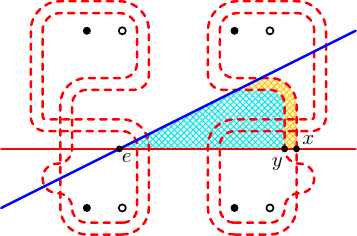}
		\caption{The cobordism maps induced by $1/2$-surgery on the disks $D$ and $D'$.}
		\label{fig:disk_surgery}
	\end{figure}
\end{example}

\begin{rem}
	For any choice of handle decompositions $\cD$ and $\cD'$ for $B^4 \smallsetminus D$, we can construct maps $F_{\cD}$ and $F_{\cD'}$ by the usual Heegaard Floer techniques.  At present, there is no general proof that the maps associated to cobordisms with corners are invariants. Though, after gluing, it follows from \cite{LOT_spectral_seq_II, zemke_linkcob} that $\bI \boxtimes F_{\cD} \simeq \bI \boxtimes F_{\cD'}$. Since we will primarily be interested in the glued maps, we will abuse notation and write $F_{B^4 \smallsetminus D}$ for any choice of map, knowing that, after pairing, the induced map is independent of the various choices in the construction of $F_{B^4 \smallsetminus D}$.

    While functoriality of the bordered Floer invariants has yet to be fully fleshed out, we expect that the maps induced by a concordance and the map computed via (any) handle decomposition of its complement specify equivalent data. Work of the second author with Sungkyung Kang \cite{guth2024invariant} considered a natural map 
    \begin{align*}
        \Lambda: \Mor^\cA(\CFDh(\KC_0), \CFDh(\KC_1)) \ra \Hom_{\cR}(\CFK_\cR(K_0), \CFK_\cR(K_1)),
    \end{align*}
    which we expect to yield a \emph{relative} Lipshitz-Ozsv\'ath-Thurston correspondence between the Heegaard Floer invariants of concordances and their complements.
\end{rem}

Recent work of Chen and Chen-Hanselman \cite{chen_hfk_satellite} has shown that the immersed curve invariants can be quite useful in understanding certain satellites. A $(1,1)$-\emph{satellite pattern} is a satellite pattern $P$ which admits a genus 1 bordered Heegaard diagram $\cH_P$. In this case, the immersed curve invariant for a knot $K$ can be paired with the beta curve of $\cH_P$ to compute a model for $\CFK^-(P(K))$. 

Given a concordance $C: K_0 \ra K_1$ and satellite pattern $P$, we can construct a \emph{satellite concordance} $P(C): P(K_0) \ra P(K_1)$ by removing a neighborhood of $C$ and replacing it with $(S^1 \times D\times I, P \times I)$. The map induced by $P(C)$ can be computed as $\bI_{P} \boxtimes F_{S^3\times I \smallsetminus K}$, according to \cite{guth_one_not_enough_exotic_surfaces}. Much as in \Cref{thm:FukayaHomPairing}, there is a nice interpretation of this in the context of immersed curves.

\begin{theorem}\label{thm:FukayaHomSatellites}
	Let $C$ be a concordance from $K_0$ to $K_1$. Let $\btheta_0$ and $\btheta_1$ be the immersed curve invariants for $\KC_0$ and $\KC_1$ respectively. Identify $F_{S^3 \times I\smallsetminus C}$ as an element $f \in \hom_{\Fuk_{\circ}}(\btheta_0, \btheta_1)$.  
	Let $P \sub S^1 \times D^2$ be a $(1,1)$-satellite pattern and let $\bbeta_P$ be its immersed curve representative. Then, the diagram
	\begin{align*}
		\begin{tikzcd}[ampersand replacement = \&]
			\hom_{\Fuk_{\circ\bullet}}(\btheta_0, \bbeta_P) \ar[r,"\simeq"] \ar[d,"m_2(f\text{,-})"] \& 
			\CFA^-(S^1\times D^2, P) \boxtimes \CFDh(\KC_0)\ar[d,"\bI_P\boxtimes F_{S^3\times I \smallsetminus C}"]\\
			\hom_{\Fuk_{\circ\bullet}}(\btheta_1, \bbeta_P)  \ar[r,"\simeq"] \& 
			\CFA^-(S^1\times D^2, P) \boxtimes \CFDh(\KC_1)
		\end{tikzcd}
	\end{align*}
	commutes up to homotopy.
\end{theorem}

\begin{example}
	In \Cref{ex:946_disks}, we determined the immersed curve invariants for the complements of the two slice disks for $m9_{46}$. Consider the \emph{Whitehead doubles} of $D$ and $D'$. By Theorem \ref{thm:FukayaHomSatellites}, this map can be computed by pairing $(\btheta_U \cup \btheta_K, \bm b)$ not with the $z$-axis representing the core circle of the solid torus, but a different curve, $\bbeta_{\Wh}$ which instead represents the Whitehead double pattern. This curve is determined by the bordered module for the Whitehead double, which has been long understood \cite{chen_hfk_satellite,levine_doubling_operators}. The new pairing diagram is shown in Figure \ref{fig:whitehead_triangle}.
	
	\begin{figure}
		\centering
		\begin{subfigure}{.4\textwidth}
			\centering
			\includegraphics[scale=1]{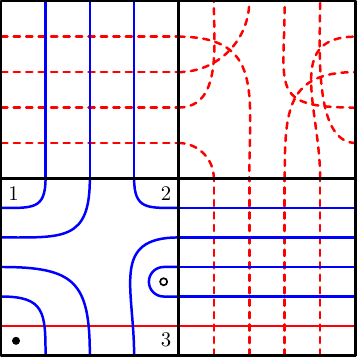}
		\end{subfigure}%
		\hspace{2cm}
		\begin{subfigure}{.4\textwidth}
			\centering
			\raisebox{-0.5cm}{\includegraphics[scale=1]{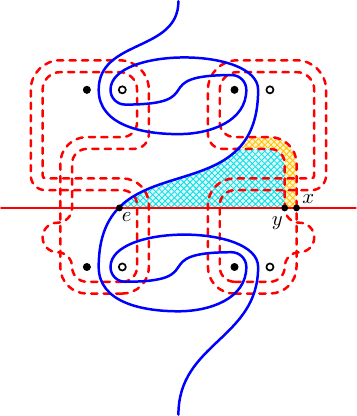}}
		\end{subfigure}
		\caption{The Whitehead double pattern (dark blue) and the two triangles specifying $F_{\Wh(D)}+F_{\Wh(D')}$ (right).}
		\label{fig:whitehead_triangle}
	\end{figure}
	
	Again, there are pairs of triangles which contribute to the map $F_{\Wh(D)} + F_{\Wh(D')}$ (see Figure \ref{fig:whitehead_triangle}); in particular, this implies that the two satellite disks are not smoothly isotopic rel boundary, recovering a result from \cite{guth2023doubled}. In \Cref{sec:applications}, we show that this phenomena is quite common, and holds for any $(1,1)$-pattern whose immersed curves satisfy a simple positivity criterion. 
\end{example}

In \cite{guth2023doubled}, Hayden--Kang--Park and the second author considered the question of the injectivity of unknotted satellite patterns: if $D$ and $D'$ are smoothly distinct slice disks for $K$, under what conditions on $P$ will $P(D)$ and $P(D')$ be smoothly distinct? There, the authors considered slice disks $D$ and $D'$ which are assumed to be distinguishable by their induced maps on $\HFKh$ and showed that $P(D)$ and $P(D')$ could be distinguished provided the bordered Floer invariants of $P$ satisfied a certain algebraic property. For $(1,1)$-satellites, this algebraic property has a simple geometric interpretation: in the universal cover, the immersed curve invariant for an unknotted pattern $P$ necessarily intersects the $x$-axis once. It turns out the $\HFK$-injectivity of $P$ only depends on the slope of this curve near this intersection with the $x$-axis. See \Cref{def:pos_slope} for a more careful exposition. 

\begin{prop}\label{prop:positive_patterns_inject}
	Let $P$ be an unknotted pattern which admits a genus 1 Heegaard diagram for which $P$ has positive slope. If $D$ and $D'$ are two $\HFKh$-distinguishable slice disks, then $P(D)$ and $P(D')$ are also $\HFKh$-distinguishable.
\end{prop}

In another direction, we consider cobordisms obtained obtained by a 4-dimensional analogue of splicing. Given a slice disk $D$ and a knot $J$, we consider the splice of $D$ along $J$, $D_J$, which we define in \Cref{sec:applications}. Much like in the satellite case, this construction can be used to construct many pairs of distinguishable 4-manifolds (including examples which are exotic rel boundary). 

\begin{prop}\label{prop:positive_tau_inject}
	Let $D$ and $D'$ be a pair of slice disks for $K$ which are distinguished by their induced map on $\HFKh$. Then, for any knot $J$ with $\tau(J) > 0$, the cobordisms $D_J$ and $D_J'$ are distinguishable by $\HFh$.
\end{prop}

We also reinterpret \cite{HL_Inv_bordered_floer,kang_bordered_involutive_HFK} in terms of immersed curves. While this reinterpretation is straightforward, at present it appears to be difficult to extract any particularly new or useful information from the immersed perspective.

\subsection*{Notation} Throughout, we will denote the field with two elements by $\F$. Additionally, since we will need to consider the space $S^1\times D^2$ regularly in later sections, we will define $\donut:=S^1\times D^2$ (for `d$\donut$nut') for the sake of notational brevity. Just as we denote the partially wrapped Fukaya category of the punctured torus with a single stop by $\Fuk_{\circ}$, we will denote the partially wrapped Fukaya category of the punctured torus with a single stop and one interior marked point by $\Fuk_{\circ\bullet}$, and the Fukaya category of the unpunctured torus with two marked points by $\Fuk_{\bullet\bullet}$. In addition, though we will encounter it only briefly, we also define $\Fuk_{\circ\circ}$ to be the partially wrapped Fukaya category of the twice-punctured torus with a single stop on one puncture.

\subsection*{Acknowledgments} The authors would like to thank Holt Bodish, Wenzhao Chen, Tobias Dyckerhoff, Matthew Habermann, Jonathan Hanselman, Jeff Hicks, Sungkyung Kang, Robert Lipshitz, Ciprian Manolescu, and Jacob Rasmussen for helpful conversations.  We would also like to specially acknowledge that Jonathan Hanselman's lecture notes for the Spring 2021 course \emph{Fukaya Categories and Floer Homology of 3-Manifolds} were very helpful during the writing of this paper.

JC gratefully acknowledges support from the Deutsche Forschungsgemeinschaft (DFG, German Research Foundation) under Germany’s Excellence Strategy - EXC 2121 ``Quantum Universe'' - 390833306. GG was supported by the Simons Collaboration Grant on New Structures in Low-Dimensional Topology. We are also grateful for the hospitality of the Max Planck Institute for Mathematics where a portion of this work was completed.

	\section{Background}
\subsection{Bordered Floer homology and immersed curves}\label{sec:background}

Let $M$ be a 3-manifold with torus boundary and choose closed, oriented, essential, geometrically dual curves $\alpha,\beta$ in $\partial M$. This pair specifies a parametrization of the boundary of $M$. 

The bordered Floer invariants of $(M,\alpha, \beta)$ are modules over the \emph{torus algebra}, $\cA$. The algebra $\mathcal{A}$ is generated over $\F$ by orthogonal idempotents $\iota_0$ and $\iota_1$, and elements $\rho_i$ for $i=1,2,3$ subject to the relations $\rho_i^2=0$ and $\rho_{i+1}\rho_i=0$ with indices taken modulo 4. 
\begin{align}
	\begin{tikzcd}[ampersand replacement=\&]
		\iota_0 \ar[rr,"\rho_1", bend left] \ar[rr,"\rho_3" above, bend right] \& \& \iota_1 \ar[ll,"\rho_2" above]
	\end{tikzcd}/ (\rho_{i+1}\rho_{i} = 0)
\end{align}
For any string $s=s_1s_2\cdots s_k$ with $s_i\in\{1,2,3\}$, we denote the product $\rho_{s_1}\rho_{s_2}\cdots\rho_{s_k}$ by $\rho_s$ whenever it is nonzero. The algebra $\mathcal{A}$ then decomposes as $\mathcal{A}=\bigoplus_{i,j=0,1}\iota_i\mathcal{A}\iota_j$, where
\begin{align}
	\begin{split}
		\iota_0\mathcal{A}\iota_0&=\langle\iota_0,\rho_{12}\rangle\\
		\iota_0\mathcal{A}\iota_1&=\langle\rho_1,\rho_3,\rho_{123}\rangle\\
		\iota_1\mathcal{A}\iota_0&=\langle\rho_2\rangle\\
		\iota_1\mathcal{A}\iota_1&=\langle\iota_1,\rho_{23}\rangle.
	\end{split}
\end{align} 

A left type-$D$ structure over $\cA$ is an $\F$-vector space $D$ with a decomposition $D = \iota_0 D \oplus \iota_1 D$ equipped with a map 
\[
\delta^1: D \ra \cA \otimes D,
\]
such the the composition 
\[
D \xra{\delta^1} \cA \otimes D \xra{\id_\cA \otimes \delta^1}\cA^{\otimes 2} \otimes D \xra{\mu_2 \otimes \id_D} D
\]
is trivial. There are also maps $\delta^k: D \ra \cA^{\otimes k} \otimes D$ defined iteratively from $\delta^1$. \footnote{One can define type-$D$ structures much more generally over an arbitrary $A_\infty$-algebra. However, the structure relation is correspondingly much more complicated --- we refer the reader to \cite[Definition 2.2.23]{LOTBimodules2015} for details.}

To a type-$D$ structure $D = \iota_0 D \oplus \iota_1 D$, there is a naturally associated directed graph. Let $\frG$ be the graph with vertices corresponding to generators of $D$ (labeled by their respective idempotents) and an edge from $x$ to $y$ labeled by $\rho_{I}$ if $\rho_I \otimes y$ appears in $\delta^1(x)$.

An $A_\infty$-module over $\cA$ is a graded $\F$-module $M$ equipped with operations 
\begin{align*}
    m_i: M \otimes \cA^{\otimes (i-1)} \ra M[2-i],
\end{align*}
satisfying 
\begin{align*}
    &\sum_{i+j = n+1} m_i(m_j(x \otimes a_1 \otimes \hdots \otimes a_{j-1}) \otimes \hdots \otimes a_{n-1}) \\
+ &\sum_{i+j=n+1}\sum_{\ell = 1}^{n-j} \mu_i(x, a_1 \otimes \hdots \otimes a_{\ell-1} \otimes \mu_j(a_\ell \otimes \hdots \otimes a_{\ell + j - 1}) \otimes\hdots \otimes a_{n-1}) = 0.
\end{align*} 
The two notions are dual in an appropriate sense (see \cite{LOT_HF_as_morphism}). This duality is manifest quite practically in the case that the type-$D$ structure is reduced (i.e. no edge in the associated train track is labeled with a $\id_\cA = \iota_0 + \iota_1)$. Namely, an oriented path from $x$ to $y$ determines a word in the alphabet $\{1,2,3\}$ by concatenating the labels of the edges making up the path and replacing every instance of a $1$ with a $3$ and vice versa. This word determines a sequence $I_1, \hdots, I_k$ by dividing $I$ into the shortest sequence of words $I_j \in \{1, 2, 3, 12, 23, 123\}.$ This data corresponds to the action 
\begin{align}
	m_{k+1}(x\otimes \rho_{I_i}\otimes \hdots \otimes \rho_{I_k})=y.
\end{align}
In particular, the graph $\frG$ encodes both $D$ and its dual type-$A$ structure $A$.

\begin{figure}[h!]
	\centering
	\includegraphics[width = .6 \linewidth]{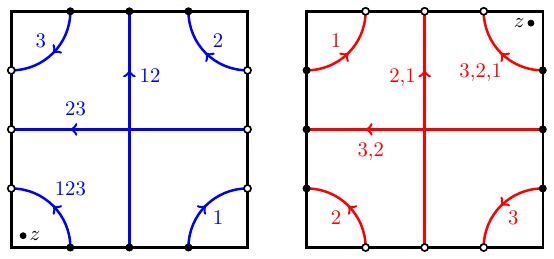}
	\caption{Conventions for constructing a train track in the torus from a type-$D$ structure.}
	\label{fig:typeD_A_conventions}
\end{figure}

The graph $\frG$ can naturally be embedded into the torus $\R^2/\Z^2$: the vertices which correspond to the generators of $\iota_0 D$ are embedded along the parametrizing curve $\alpha$, which we identify with the segment $y = 0$ in $\R^2/\Z^2$, and those corresponding to the generators of $\iota_1 D$ along $\beta$, which we identify with $x = 0$; the edges are embedded in $[0,1]\times[0,1]$ according to their labelings by algebra elements as shown in Figure \ref{fig:typeD_A_conventions}.
Following \cite{hanselman2023bordered}, we refer to this embedding as the \emph{type-$D$ realization} of $\frG$, or $D(\frG).$ The \emph{type-$A$ realization} is obtained by reflecting $D(\frG)$ along the line $y = -x$ and swapping $\rho_1$ and $\rho_3$. 

\begin{figure}
	\centering
	\includegraphics[scale=1]{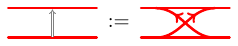}
	\caption{A crossover arrow.}
	\label{fig:crossover_arrow}
\end{figure}

\begin{figure}
	\centering
	\includegraphics[width=.8 \linewidth]{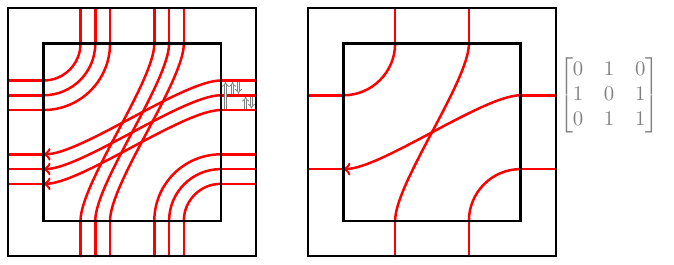}
	\caption{A local system presented as a collection of crossover arrows.}
	\label{fig:crossover_localsys}
\end{figure}

A key result of \cite{hanselman2023bordered} is that for a type-$D$ structure associated to a 3-manifold with torus boundary, the associated immersed train track can be assumed to take a particularly nice form: namely, there is a representative of the homotopy equivalence type of the type-$D$ structure for which the associated train track is an immersed curve decorated by local systems, which can be encoded by a collection of \emph{cross-over arrows} (see Figures \ref{fig:crossover_arrow} and \ref{fig:crossover_localsys}). The curve invariant for $\CFDh(M)$ is defined by representing the module as such a graph and including it into $\partial M$ according to the conventions above. 

A type $A$ structure $A$ can be paired with a type $D$ structure, $D$. Define $A \boxtimes D$ to be the chain complex $A \otimes_{\cI} D$ equipped with the differential 
\[
\partial^\boxtimes(x\otimes y) = \sum_{k=0}^\infty (m_{k+1}\otimes \id_D)(x\otimes \delta^k(y)).
\]
If $\frG_1$ and $\frG_2$ are $\cA$-decorated train tracks representing a type-$A$ structure $A$ and type-$D$ structure $D$, then the Lagrangian intersection Floer homology of $A(\frG_1)$ and $D(\frG_2))$ exactly recovers $A \boxtimes D$. By construction, there is a bijection between the two generating sets, and terms in the differential $\partial^\boxtimes$ correspond precisely to bigons connecting intersection points between $A(\frG_1)$ and $D(\frG_2)$.

If $M$ and $N$ are 3-manifolds with torus boundary whose bordered invariants are represented by immersed curves $\btheta_M$ and $\btheta_N$ and $h: \partial N \ra \partial M$ is an orientation reversing diffeomorphism, Hanselman--Rasmussen--Watson use this construction to define a pairing $\langle \btheta_M, \btheta_N \rangle$ which recovers the Floer homology of $M \cup_h N$ as follows. We fix a parametrization $(\alpha', \beta')$ for $\partial N$, and parametrize $\partial M$ by $(\alpha = h(\beta'), \beta = h(\alpha'))$. We may therefore view $\btheta_M$ and $\btheta_N$ as residing in the same parametrized torus. However, with respect to our chosen parametrizations, the map $h$ swaps $\alpha$ with $\beta'$ and $\beta$ with $\alpha'$, which amounts to reflecting $\btheta_N$ across the line $y = x$. This is almost the type D realization, which is given by reflection across the line $y = -x$; to utilize the construction above, we must compose $h$ with the elliptic involution. This composition is denoted $\overline{h}$. In other words, the Lagrangian Floer homology of $\btheta_M$ with $\overline{h}(\btheta_N)$ computes exactly $\CFAh(M) \boxtimes \CFDh(N) \simeq \CFh(M \cup_h N).$ In \cite{hrw_properties_apps}, it is shown that the elliptic involution acts by the identity on unlabeled curves, so this subtlety is largely cosmetic. 

We shall take a slightly different perspective. We will always identity $\partial M$ and $\partial N$ via an orientation \emph{preserving} diffeomorphism $h: \partial N \ra \partial M$, in order to recover the morphism pairing theorem, which states that 
\begin{align*}
    \Mor_\cA(\CFAh(M), \CFAh(N)) \simeq \CFh(-M \cup_h N).
\end{align*}
In the setting of immersed curves, this corresponds to fixing a parametrization $(\alpha', \beta')$ for $\partial N$, and parametrizing $\partial M$ and $(\alpha = h(\alpha'), \beta = h(\beta'))$. As we shall see, $\hom_{\Fuk_\circ}(\btheta_M, h(\btheta_N))$ does indeed compute $\Mor_\cA(\CFAh(M), \CFAh(N))$, and, curiously, the elliptic involution does not appear.

\subsection{The partially wrapped Fukaya category of the punctured torus}
We recall here a cursory definition of the partially wrapped Fukaya category, omitting most technical details. Suppose that $(M,\omega)$ is a symplectic manifold with contact boundary and let $\hat{M}=M\cup_{\partial M}[1,\infty)\times\partial M$ be the completion of $M$ by the positive symplectization $[1,\infty)\times\partial M$ of the boundary. Suppose that $\bsigma\subset\partial M$ is a collection of hypersurfaces (by abuse of notation, we will identify $\bsigma$ with $\bigsqcup_{\sigma\in\bsigma}\sigma$). Choose a smooth function $\rho:\partial M\to[0,1]$ such that $\bsigma=\rho^{-1}(0)$ and let $H_\rho:\hat{M}\to\R$ be a non-negatively valued Hamiltonian given on the symplectization $[1,\infty)_r\times\partial M_x$ by $H_\rho(r,x)=\rho(x)r$.
\begin{definition}
	The \emph{partially wrapped Fukaya category} $\Fuk_{\bsigma}(M)$ is the $A_\infty$-category whose objects are (possibly compact, immersed) oriented Lagrangian submanifolds of $\hat{M}$ with cylindrical ends of the form $[1,\infty)\times\Lambda$ for $\Lambda\subset\partial M\smallsetminus\bsigma$ a Legendrian submanifold. Morphism spaces in $\Fuk_{\bsigma}(M)$ are defined by
	\begin{align}
		\hom_{\Fuk_{\bsigma}}(L_1,L_2)=\lim\limits_{w\to+\infty}\CF(\phi_{wH_\rho}(L_1),L_2),
	\end{align}
	where $\phi_{wH_\rho}$ is the Hamiltonian flow of $wH_\rho$, with differential, composition, and higher composition maps defined by counting $H_\rho$-perturbed holomorphic polygons.
\end{definition}
\begin{remark}
	In our setting, with $M=\Sigma$ a surface, we also equip our Lagrangians with $\F$-local systems, which we may equivalently think of as $\F$-vector bundles. If $\bm{L}_i=(L_i,E^i)$ are pairs consisting of a Lagrangian $L_i$ together with a local system $E^i$, the definition of the Lagrangian Floer complex must be modified as follows. Here,
	\begin{align}
		\CF(\bm{L}_0,\bm{L}_1)=\bigoplus_{x\in L_1\cap L_2}\Hom_{\F}(E^0_x,E^1_x)
	\end{align}
	and the differential is defined on $f\in\Hom_{\F}(E^0_x,E^1_x)$ by
	\begin{align}
		\partial f=\sum_{y\in L_1\cap L_2}\sum_{u\in\pi_2(x,y),\,\mu(u)=1}P_1\circ f\circ P_0,
	\end{align}
	where, for a homotopy class $u$ of strips from $x$ to $y$ with a holomorphic representative $\R\times[0,1]\to\Sigma$, $P_i$ is the parallel transport map of the local system on $L_i$ along the path $\gamma_i^u:[-\infty,\infty]\to L_i$ given by $\gamma_i^u(t)=u((-1)^it,i)$. In other words, $\gamma_0$ is the path in $L_0$ from $x$ to $y$ given by the boundary of $u$ in $L_0$ and $\gamma_1$ is the path in $L_1$ from $y$ to $x$ given similarly.
	\begin{figure}
		\centering
		\includegraphics[width=\textwidth]{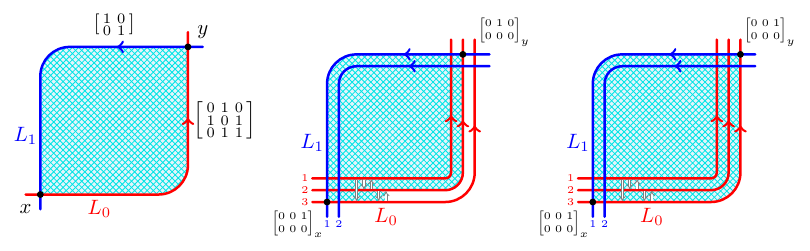}
		\caption{Disks between Lagrangians $L_0$ and $L_1$ with local systems (with parallel transport on each segment represented by matrices) of rank 3 and 2, respectively, and one of the corresponding disks between their train track representatives which shows that $\partial\left[\begin{smallmatrix}
				0 & 0 & 1\\
				0 & 0 & 0
			\end{smallmatrix}\right]_x=\left[\begin{smallmatrix}
			0 & 1 & 1\\
			0 & 0 & 0
		\end{smallmatrix}\right]_y+(\textup{other terms})$ has $\left[\begin{smallmatrix}
			0 & 1 & 0\\
			0 & 0 & 0
		\end{smallmatrix}\right]_y$ and $\left[\begin{smallmatrix}
		0 & 1 & 1\\
		0 & 0 & 0
	\end{smallmatrix}\right]_y$ as summands.}\label{fig:LocalSystemDisk}
	\end{figure}
	Since we are working with curves and $\F$-local systems, we may instead think of pairs $(L,E)$ as an immersed train track by replacing $L$ with $\mathrm{rk}E$ parallel copies of $L$ and placing crossover arrows between parallel strands according to the parallel transport of $E$ along $L$, given some choice of orientation on $L$, as in \cite[Section 3]{hanselman2023bordered}. When we have two such pairs $(L_0,E^0)$ and $(L_1,E^1)$, we can arrange so that these crossover arrows are placed on parallel segments corresponding to segments lying between intersection points according to the parallel transport maps. Numbering the parallel strands by $1,\dots,\mathrm{rk}E^i$ allows us to identify the span of the $(\mathrm{rk}E^0)(\mathrm{rk}E^1)$ intersection points between these immersed train tracks corresponding to an intersection point $x\in L_0\cap L_1$ with the space of matrices $\mathrm{Mat}_{\mathrm{rk}E^1\times\mathrm{rk}E^0}(\F)\cong\Hom_{\F}(E^0_x,E^1_x)$. If $x$ and $y$ are connected by a holomorphic strip and $f\in\Hom_{\F}(E^0_x,E^1_x)$ is given by an elementary matrix, the strips with boundary on the train track representatives which emanate from the intersection point corresponding to $f$, end on an intersection point coming from $y$, and whose boundaries are consistently oriented paths with respect to the chosen orientations of the Lagrangians and crossover arrows, then give us the summands of $\partial f$ under this identification. See Figure \ref{fig:LocalSystemDisk} for an example. In the case that both local systems are of rank 1 --- which will typically be the case in our setting --- this construction recovers the ordinary Lagrangian Floer complex. For a more thorough account of $\F$-local systems in Lagrangian Floer homology, we refer the reader to \cite{Konstantinov}.
\end{remark}
To give a complete definition of partially wrapped Fukaya categories with immersed objects in full detail is beyond the scope of this paper. We instead refer the reader to \cite[Appendix A]{auroux_fukaya_sym_products} for further details on the full subcategory generated by embedded Lagrangians in the exact case. However, we will be concerned primarily with the case that $M$ is the surface $\Sigma=T^2\smallsetminus D$, where $D$ is a small open disk, and $\bsigma=z\in\partial \Sigma$ is a single point.\footnote{Occasionally, we will also need to consider the same surface with an additional puncture but no additional stops, or with a marked point in the interior, and sometimes just $T^2$ with two marked points} In this case, we will denote the partially wrapped Fukaya category simply by $\Fuk_{\circ}$. In this setting, the Lagrangians $\phi_{wH_\rho}(L_1)$ and $L_2$ become transverse and no new intersection points are created for $w$ sufficiently large so we may interpret the hom complex $\hom_{\Fuk_{\circ}}(L_1,L_2)$ as an ordinary Lagrangian Floer complex $\CF(\phi(L_1),L_2)$, where $\phi(L_1)$ is the result of dragging any endpoints of $L_1$ which lie on $\partial\Sigma$ around the boundary according to the orientation of the surface until they are near $z$ in such a way that further wrapping will not introduce new intersection points.
\begin{figure}
	\begin{center}
		\includegraphics[width=\textwidth]{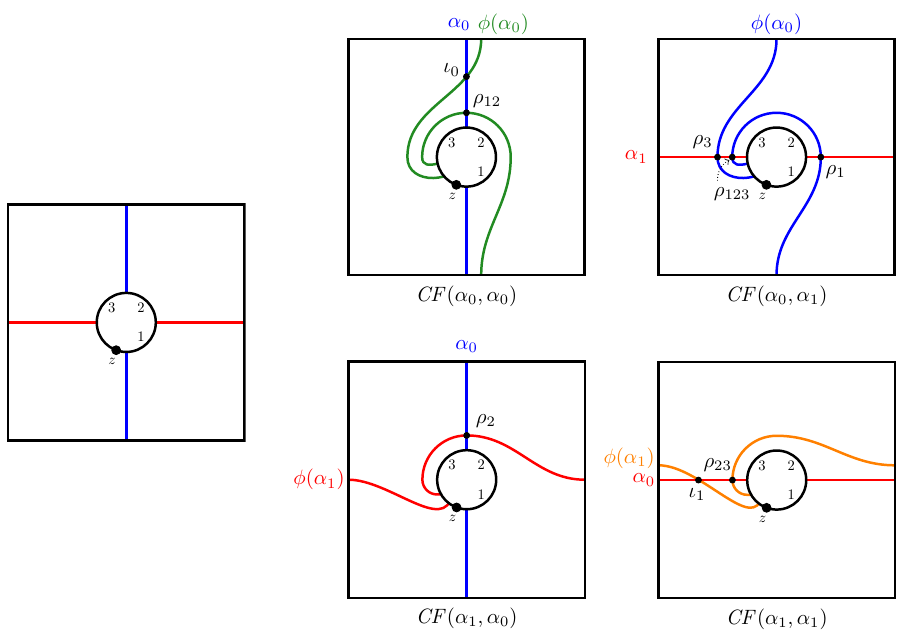}
	\end{center}
	\caption{The generators of $\mathcal{A}$ as morphisms in the partially wrapped Fukaya category $\Fuk_{\circ}$.}
	\label{fig:A-as-Fukaya}
\end{figure}
Fix a handle decomposition of the torus consisting of a single 0- and 2-handle and two 1-handles and identify the disk $D$ as above with the interior of the 0-handle in this decomposition. The cores of the 1-handles determine two arcs, $\alpha_0$ and $\alpha_1$, in $\Sigma$. By work of Auroux \cite{auroux_fukaya_sym_products}, $\Fuk_{\circ}$ is split-generated by the arcs $\alpha_0$ and $\alpha_1$, so Lagrangians in $\Sigma$ may be viewed as modules over the endomorphism $A_\infty$-algebra
\begin{align}
	\End_{\Fuk_{\circ}}(\alpha_0\oplus\alpha_1)=\bigoplus_{i,j=0,1}\hom_{\Fuk_{\circ}}(\alpha_i,\alpha_j).
\end{align}
Moreover, Auroux shows that $\End_{\Fuk_{\circ}}(\alpha_0\oplus\alpha_1)$ is isomorphic as an $A_\infty$-algebra, with differential and composition maps given by counting holomorphic polygons with Lagrangian boundary conditions, to the algebra $\mathcal{A}$ associated to the torus by bordered Floer homology (see Figures \ref{fig:A-as-Fukaya} and \ref{fig:TriangleProduct}). The identification of each intersection point with a generator of $\mathcal{A}$ proceeds as follows: if, in the process of wrapping, the endpoint of $\phi(\alpha_i)$ traverses a sequence of Reeb chords $\rho_{i_1},\dots,\rho_{i_k}$ in $\partial\Sigma$ --- as indicated by the numbering of the boundary circle in Figure \ref{fig:A-as-Fukaya} --- before crossing $\alpha_j$, then the intersection point introduced when the endpoint crosses the endpoint of $\alpha_j$ is identified (uniquely) with the product $\rho_{i_1}\cdots\rho_{i_k}$. Otherwise, the intersection point under consideration was introduced as the first intersection of $\phi(\alpha_i)$ with $\alpha_i$ and is identified with $\iota_i$. Moreover, this identification makes the left- and right-idempotents of the algebra elements $\rho_s$ clear: as is to be expected, if $\rho_s$ is identified with an intersection of $\phi(\alpha_i)$ with $\alpha_i$, then $\rho_s\in\iota_i\mathcal{A}\iota_j$.
\begin{figure}
	\begin{center}
		\includegraphics[]{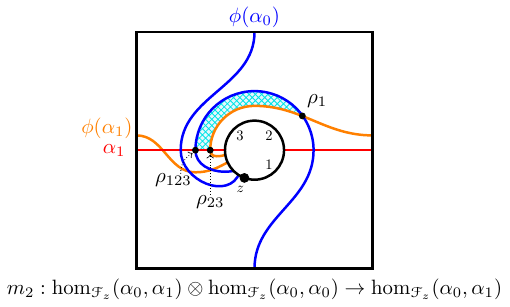}
	\end{center}
	\caption{A holomorphic triangle furnishing the product $m_2(\rho_1\otimes\rho_{23})=\rho_{123}$.}
	\label{fig:TriangleProduct}
\end{figure}
By \cite[Proposition 6.8]{auroux_fukaya_sym_products}, the Yoneda embedding $\yo:\Fuk_{\circ}\to\mathrm{mod-}\mathcal{A}$ given by
\begin{align}
	L\mapsto\yo(L):=\hom_{\Fuk_{\circ}}(L,\alpha_0\oplus\alpha_1)
\end{align}
is a cohomologically full and faithful contravariant embedding\footnote{Since the punctured torus has genus 1, there is a canonical identification of $\Fuk_{\circ}$ with the extended Fukaya category $\Fuk_{\circ}^\sharp$ of Lagrangian correspondences used in \cite{auroux_fukaya_sym_products}.}.

The work of \cite{hanselman2023bordered} shows that the $A_\infty$ category $tw\Fuk_{\circ}$ of twisted complexes over $\Fuk_{\circ}$ is equivalent to the Fukaya category of immersed curves in the punctured torus with local systems. Since the curves $\alpha_0$ and $\alpha_1$ split-generate $\Fuk_{\circ}$, this identifies the latter Fukaya category with the category of right $A_\infty$-modules over $\End_{\Fuk_{\circ}}(\alpha_0\oplus\alpha_1)$.

	\section{Immersed Curves and Composition}
\subsection{Composition for 3-manifolds}

The composition theorem for 3-manifolds with torus boundary is a straightforward extension of the work of Hanselman--Rasmussen--Watson \cite{hanselman2023bordered}. Recall that the immersed curve invariant $\btheta_M$ for $M$ is \emph{by definition} a model for $\CFAh(M)$, represented as a graph which has been embedded into the punctured torus. On the other hand, via the Yoneda embedding functor, we can also pair $\btheta_M$ with the arcs $\alpha_0$ and $\alpha_1$; the resulting object $\yo(\btheta_M)$ is an $A_\infty$-module over $\End_{\Fuk_\circ}(\alpha_0\oplus \alpha_1)$, which by the work of Aurox, is identified with the torus algebra. Unsurprisingly, these two modules are equivalent. 

\begin{proposition}\label{prop:CFA is Y(gamma)}
	Let $M$ be a bordered 3-manifold with torus boundary and let $\btheta_M$ be its decorated immersed multi-curve invariant. Then, there is a homotopy equivalence 
	\begin{align}
		\yo(\btheta_M)\simeq\CFAh(M)
	\end{align}
	of right $A_\infty$-modules.
\end{proposition}

\begin{proof}
	This is ultimately a reformulation of the work in \cite{hanselman2023bordered} and closely follows the ideas in their work. 

    Recall that the curve $\btheta_M$ is a graph representing $\CFAh(M)$ whose vertices form a basis for $\CFAh(M)$; this graph is embedded into the solid torus such that the vertices are embedded along the parametrizing curves $\alpha$ and $\beta$. The generators of $\yo(\btheta_M)$ are precisely the intersection points between $\btheta_M$ and the two alpha arcs $\alpha_0$ and $\alpha_1$. The generators of these two modules are clearly identified.

    Let $\mu_k$ and $m_k$ denote the $A_\infty$-operations for $\yo(\btheta)$ and $\CFAh(M)$ respectively. It suffices to identify the actions of the two modules, i.e. to show that terms in $m_{k+1}:\CFAh(M)\otimes\cA^{\otimes k}\ra\CFAh(M)$
	correspond uniquely to $(k+2)$-gons contributing to the map $\mu_{k+1}:\yo(\btheta) \otimes\End_{\Fuk_{\circ}}(\alpha_0 \oplus \alpha_1)^{\otimes k}\ra\yo(\btheta)$
	under the identification of $\cA$ with $\End_{\Fuk_{\circ}}(\alpha_0\oplus\alpha_1).$ The $A_\infty$-operations on $\CFAh(M)$ can be read off from paths in $\btheta_M$. Let $x$ and $y$ be generators of $\CFAh(M)$ (i.e. intersection points between $\btheta_M$ and $\alpha_0 \oplus \alpha_1$.) An oriented path from $x$ to $y$ determines a word in the alphabet $\{1,2,3\}$ by concatenating the labels of the edges making up the path. This word determines a sequence $I_1, \hdots, I_k$ by dividing $I$ into the shortest sequence of words $I_j \in \{1, 2, 3, 12, 23, 123\}.$ This data corresponds to the action 
	\begin{align}
		m_{k+1}(x\otimes \rho_{I_i}\otimes \hdots \otimes \rho_{I_k})=y.
	\end{align}    
	The $A_\infty$-operations of $\yo(\btheta_M)$ are given by counting immersed $k$-gons.
    
	Just as in \cite[Theorem 2.2]{hanselman2023bordered}, any immersed polygon in the punctured torus can be decomposed as a union of embedded pieces. Since these polygons cannot pass the basepoint, there is a very restricted set of basic building blocks. Moreover, any immersed polygon which arises as a term in the action of $\End_{\Fuk_{\circ}}(\alpha_0\oplus \alpha_1)$ has exactly two vertices which appear as intersections between $\btheta$ and the $\alpha$ arcs, and the remaining $(k-1)$ vertices appear close to the puncture, labeled by algebra elements.
	
	First, suppose there is a $(k+2)$-gon from $x$ to $y$. This polygon defines two paths from $x$ to $y$, one along $\btheta$ and the other along the $\alpha$-arcs and the (wrapped) $\phi(\alpha)$-arcs (see Figure \ref{fig:IntersectionWrapped}).
	
	\begin{figure}
		\begin{center}
			\includegraphics{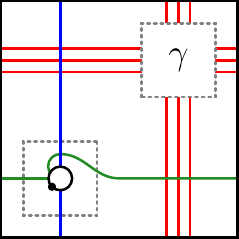}
		\end{center}
		\caption{Intersections between $\btheta$ and the (wrapped) $\alpha$-arcs.}
		\label{fig:IntersectionWrapped}
	\end{figure}
	
	As described above, the path in $\btheta$ records a sequence, $J_1,J_2,\dots,J_N$, where each $J_i$ is an element of $\{1, 2, 3, 12, 23, 123\}$, corresponding to an operation
	\begin{align}
		m_{N+1}(x \otimes \rho_{J_1} \otimes \hdots \otimes \rho_{J_N}) = y.
	\end{align}
	The path in the $\alpha$ arcs also records a sequence $I_1,\dots,I_{k}$, where the $j$-th vertex is labeled by an algebra element $\rho_{I_j}$, giving rise to an action 
	\begin{align}
		\mu_{k+1}(x\otimes\rho_{I_1}\otimes\cdots\otimes \rho_{I_k}) = y.
	\end{align}
	However, just as in the proof in \cite[Theorem 2.2]{hanselman2023bordered} (by considering the possible decompositions of polygons in Figure \ref{fig:PolygonDecompositions}) we see that $I_j=J_j$ for all $j$ and $N=k$. Therefore, a $(k+1)$-gon witnesses a $m_{k+1}$ operation on $\CFAh(M)$.
	\begin{figure}
		\begin{center}
			\includegraphics[width=\textwidth]{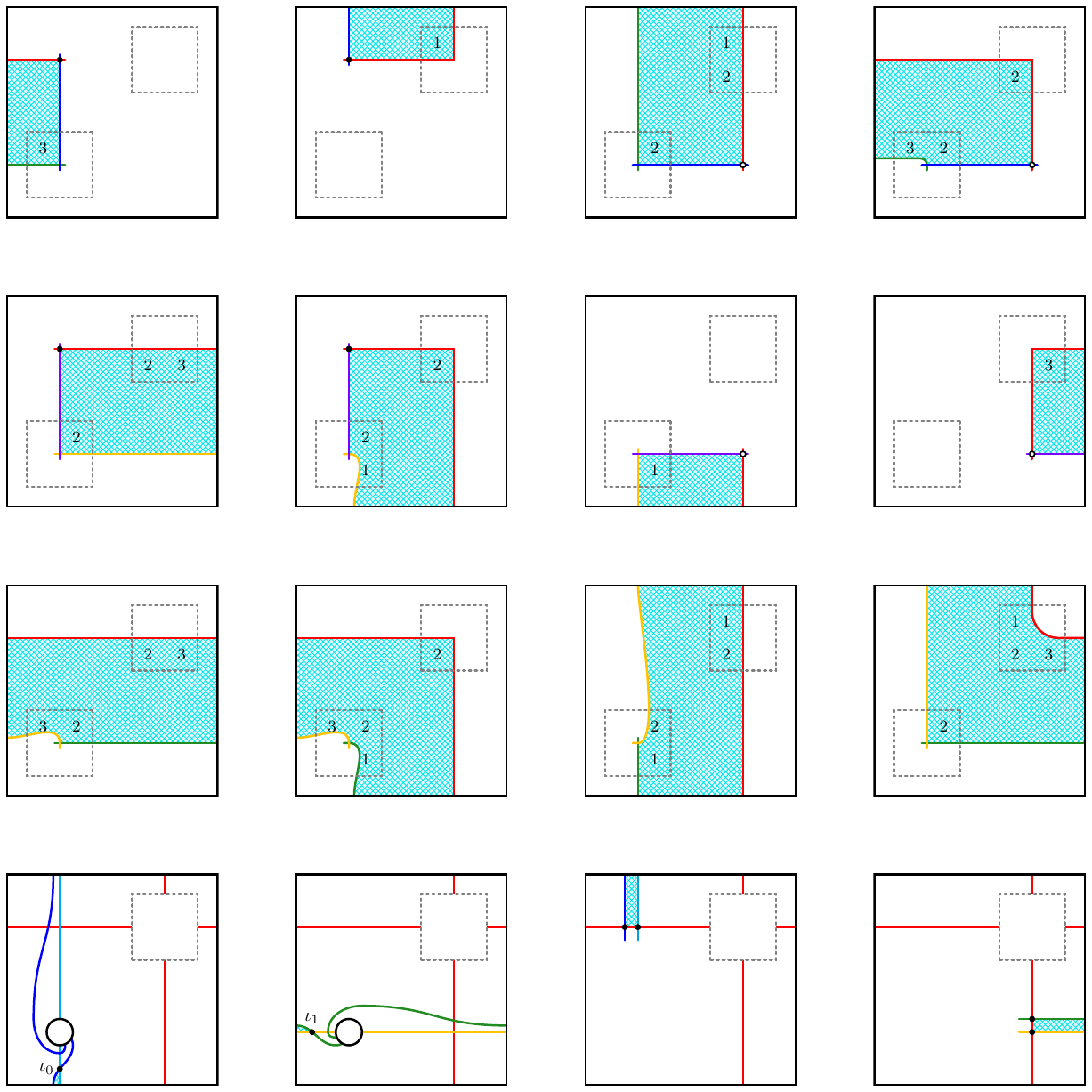}
		\end{center}
		\caption{Decompositions of polygons in $\Sigma$.}\label{fig:PolygonDecompositions}
	\end{figure}
	Conversely, if $m_{k+1}(x\otimes\rho_{I_1}\otimes\cdots\otimes\rho_{I_k}) = y$ in $\CFAh(M)$, then, by construction of the immersed curve invariant, there is a path in $\btheta$ from $x$ to $y$ which records the string $I_1\cdots I_{k-1}.$ This string determines a piece-wise linear path along the $\alpha$ arcs from $x$ to $y$ with corners on the intersection points labeled by $I_1,\dots,I_{k-1}.$ These two paths form the boundary of an immersed $(k+2)$-gon from $x$ to $y$, contributing a term $\mu_{k+1}(x\otimes\rho_{I_1}\otimes\cdots\otimes\rho_{I_k})=y$ in $\mu_k$. 
	
	Therefore, our chosen models for $\yo(\btheta)$ and $\CFAh(M)$ are genuinely isomorphic as $A_\infty$-modules. It then follows that any two models for these objects will be homotopy equivalent in general. In the presence of higher rank local systems, the same argument applies after passing to immersed train track representatives.
\end{proof}

The proof of Theorem \ref{thm:FukayaHomPairing} now follows from the Yoneda embedding construction. 

\begin{proof}[Proof of Theorem \ref{thm:FukayaHomPairing}:]
	Since $\yo$ is a fully faithful embedding and $\cA\cong\End_{\Fuk_{\circ}}(\alpha_0\oplus\alpha_1)$, Propsition \ref{prop:CFA is Y(gamma)} yields
	\begin{align}
		\begin{split}
			\Mor_\cA(\CFAh(M_1),\CFAh(M_2)) & \simeq \Mor_{\End_{\Fuk_{\circ}}(\alpha_0\oplus\alpha_1)}(\yo(\btheta_1),\yo(h(\btheta_2)))\\
			& \simeq \hom_{\Fuk_{\circ}}(\btheta_1, h(\btheta_2)).
		\end{split}
	\end{align}
	Functoriality of $\yo$ then implies the desired result. 
\end{proof}

The identification of the two hom-spaces themselves is not new, and follows from (the proof of) the bordered Floer homology morphism pairing theorem, \cite[Theorem 1]{LOT_HF_as_morphism}, and the pairing theorem for curves \cite[Theorem 1.2]{hanselman2023bordered}, which tell us that if $h: \partial N \ra \partial M$ is an orientation preserving diffeomorphism, then
\begin{align}
	\Mor_\cA(\CFAh(M), \CFAh(N)) \simeq \CFAh(-M)\boxtimes  \CFDh(N) \simeq \hom_{\Fuk_{\circ}}(\btheta_{-M}, h(\btheta_N)).
\end{align}
However, the proof above gives an \emph{explicit} identification of the two spaces: an intersection point $p$ between $\btheta_M$ and $h(\btheta_N)$ determines a morphism
\begin{align*}
    \CFAh(M) \simeq \yo(\theta_M) \ra \yo(h(\btheta_N)) \simeq \CFAh(N)
\end{align*}
by counting triangles with a corner at $p$. See \Cref{ex:MorphismTrianglesExample} for an example of this. Moreover, the proof of \Cref{thm:FukayaHomPairing} identifies  the composition in the bordered setting with composition in $\Fuk_{\circ}$.

Further, as a consequence, the analogous statement holds for type-$D$ structures.
\begin{corollary}\label{cor:pairing type D}
	If $M_i$ are 3-manifolds with torus boundary with corresponding type-$D$ immersed curve invariants $\btheta_i$, then $\CFDh(M_i)\simeq\yo^{\mathrm{op}}(\btheta_i):=\hom_{\Fuk_{\circ}}(\alpha_0\oplus\alpha_1,\btheta_i)$ and there are homotopy equivalences of chain complexes
	\begin{align}
		\Mor^\cA(M_i,M_j)\simeq\hom_{\Fuk_{\circ}}(\btheta_i,\btheta_j)
	\end{align}
	such that the squares
	\begin{align}
		\begin{tikzcd}[ampersand replacement=\&]
			\Mor^\cA(M_0,M_1)\otimes\cdots\otimes\Mor^\cA(M_{k-1},M_k)\arrow[r,"\circ_k"]\arrow[d,"\simeq"'] \& \Mor^\cA(M_0,M_k)\arrow[d,"\simeq"]\\
			\hom_{\Fuk_{\circ}}(\btheta_0,\btheta_1)\otimes\cdots\otimes\hom_{\Fuk_{\circ}}(\btheta_{k-1},\btheta_k)\arrow[r,"m_k"] \& \hom_{\Fuk_{\circ}}(\btheta_0,\btheta_k)
		\end{tikzcd}
	\end{align}
	commute up to homotopy. Note: for $k>2$, the higher composition maps $\circ_k$ in the top row all vanish because $\cA$ is a genuine dg-algebra.
\end{corollary}
\begin{proof}
	This follows from the type-$A$ case since the reflection across the line $y=-x$ which recovers a type-$D$ curve invariant $\btheta$ from its type-$A$ counterpart $\btheta^\vee$ furnishes a contravariant autoequivalence $\Fuk_{\circ}\to\Fuk_{\circ}$ which induces an equivalence
	\begin{align}
		\Mor_{\End_{\Fuk_{\circ}}(\alpha_0\oplus\alpha_1)}(\yo(\btheta_1^\vee),\yo(\btheta_2^\vee))\simeq\Mor^{\End_{\Fuk_{\circ}}(\alpha_0\oplus\alpha_1)}(\yo^{\mathrm{op}}(\btheta_1),\yo^{\mathrm{op}}(\btheta_2))
	\end{align}
	and a similar argument to that given in Proposition \ref{prop:CFA is Y(gamma)} identifies $\yo^{\mathrm{op}}(\btheta_i)$ with $\CFDh(M_i)$.
\end{proof}
This also immediately implies the following corollary.
\begin{corollary}\label{cor:MorphismIdentification}
	Let $\btheta_1$ and $\btheta_2$ be immersed curves representing type-$D$ structures $D_1$ and $D_2$. Let $f$ be an intersection point of $\btheta_1$ and $\btheta_2$. Then, under the identification above, $f$ corresponds to the module homomorphism 
	\begin{align}
		m_2(-,f):\yo^{\mathrm{op}}(\btheta_1)\to\yo^{\mathrm{op}}(\btheta_2),
	\end{align}
	regarded as an element of $\Mor^\cA(D_1,D_2)$.\hfill\qedsymbol
\end{corollary}
\begin{remark}
	We will frequently conflate the maps $\circ_2(-,f)$ and $m_2(-,f)$. We will also often abuse notation and refer to the module homomorphisms simply by $f$.
\end{remark}

Consequently, for computations, we will often work instead with morphisms of type D structures rather than morphisms of type A structures. 

\begin{figure}
	\centering
	\includegraphics[width=0.325\textwidth]{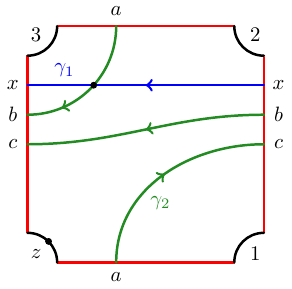}\hspace{1cm}\includegraphics[width=0.325\textwidth]{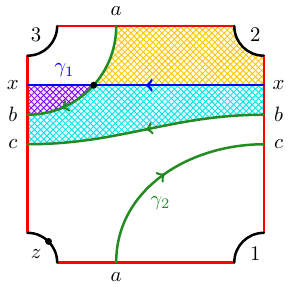}
	\caption{Immersed curves $\btheta_1$ and $\btheta_2$ representing the type-$D$ structures $D_1$ and $D_2$ in Example \ref{ex:MorphismTrianglesExample}, and the three holomorphic triangles with corners at the unique intersection point $f\in\btheta_1\cap\btheta_2$. Note: the arrows on the segments of the $\btheta_i$ correspond to the type-$D$ structure map and do not represent orientations of the curves.}
	\label{fig:int_as_morph}
\end{figure}

\begin{example}\label{ex:MorphismTrianglesExample}
	As a simple example, consider the immersed curves $\btheta_1$ and $\btheta_2$ in Figure \ref{fig:int_as_morph} representing the type-$D$ structures
	\begin{align*}
		D_1=\begin{tikzcd}[ampersand replacement=\&]
			x\arrow[loop right,out=35,in=-35,looseness=10,"23"]
		\end{tikzcd}
		\quad  \quad
		D_2=\begin{tikzcd}[ampersand replacement=\&,column sep=0.25cm]
			a\arrow[rr,"3"]\arrow[dr,"1"'] \& \& b\arrow[dl,"23"]\\
			\& c \&
		\end{tikzcd}
	\end{align*}
	over $\cA$, respectively.
	There is a single intersection point, $f$, between $\btheta_1$ and $\btheta_2$ --- which is, therefore, necessarily a cycle --- and there are three generalized holomorphic triangles with a corner at $f$ (see Figure \ref{fig:int_as_morph}.). Here, rather than draw the wrapped curves, we draw in the Reeb chords, as in \cite{LOT_bordered_HF}. The corresponding summands of the cycle $f\in\Mor^\cA(D_1,D_2)$ are
	\begin{align*}
		\raisebox{-0.5cm}{\includegraphics[]{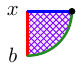}}:x\mapsto b,
	\end{align*}
	\begin{align*}
		\raisebox{-0.625cm}{\includegraphics[]{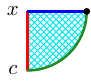}}:x\mapsto c,
	\end{align*}
	and
	\begin{align*}
		\raisebox{-0.625cm}{\includegraphics[]{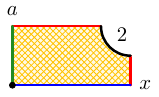}}:x\mapsto\rho_2a.
	\end{align*}
	Indeed, we have $H_*\Mor^\cA(D_1,D_2) \cong \F\langle [x\mapsto b+c+\rho_2 a]\rangle$ by direct computation.
\end{example}

	\section{Composition and Knot Floer Homology}
\subsection{Composition for satellites}
We turn to the proof of Theorem \ref{thm:FukayaHomSatellites}. Throughout the section, we will write $(T_\infty, \nu)$ for the standard doubly pointed bordered Heegaard diagram representing an $S^1$-fiber in the the $\infty$-framed solid torus. The bordered Floer homology pairing theorem \cite[Theorem 11.19]{LOT_bordered_HF} states
\begin{align}
	\CFA^-(T_\infty, \nu)\boxtimes \CFDh(\KC)& \simeq \CFK^-(S^3, K)
\end{align}

More generally, if $P$ is a pattern knot in the solid torus represented by a doubly pointed bordered Heegaard diagram $\cH_P$, the pairing theorem implies that
\begin{align}
	\CFA^-(\cH_P)\boxtimes \CFDh(\KC) \simeq \CFK^-(S^3, P(K)).
\end{align}
In the case that the pattern $P$ admits a genus 1 bordered Heegaard diagram, this pairing has an immersed curve interpretation due to Chen \cite{chen_hfk_satellite}. A genus 1 bordered Heegaard diagram for $P$ consists of two $\alpha$ arcs and a single $\beta$ circle (which will be denoted $\bbeta_P$). By construction, $\bbeta_P$ is an immersed (embedded, in fact) Lagrangian in the punctured torus with a marked point. Therefore, $\bbeta_P$ can naturally be paired with an immersed multicurve. \cite[Theorem 1.2]{chen_hfk_satellite} states that when $\btheta_K$ is the immersed curve invariant for $\KC$, then
\begin{align*}
	\hom_{\Fuk_{\circ\bullet}}(\btheta_K, \bbeta_P) \simeq \CFK^-(S^3, P(K)),
\end{align*}
where $\Fuk_{\circ\bullet}$ is the Fukaya category of the once-punctured torus with one marked point. As the notation indicates, Chen's original work only recovers $\CFK^-$ from the immersed curves. However, this result was improved in \cite{chen_hanselman_2023satellite} to show that 
\begin{align*}
	\hom_{\Fuk_{\bullet\bullet}}(\btheta_K, \bbeta_P) \simeq \CFK_\cR(S^3, P(K)),
\end{align*}where $\Fuk_{\bullet\bullet}$ is the Fukaya category of the unpunctured torus with two marked points, by working with \emph{immersed Heegaard diagrams}.

Of course, there ought to be four-dimensional analogues of these results. Given a concordance, $C: K_0 \ra K_1$ be a concordance, there is an induced map 
\[
F_{S^3\times I \smallsetminus C}: \CFDh(\KC_0) \ra \CFDh(\KC_1).
\]
Under the pairing theorem, the map induced by satellites of $C$ can be computed as 
\[
\bI_P \boxtimes F_{S^3\times I \smallsetminus C} \simeq F_{P(C)}.
\]
See \cite{guth_one_not_enough_exotic_surfaces} for a proof.
Naturally, one would expect this map to be computable in terms of immersed curves. 

It will turn out to be useful to combine the morphism and tensor product pairing theorems. Rather than view composition as a map 
\begin{align*}
    \Mor^\cA(M_0, M_1) \otimes \Mor^\cA(M_1, M_2) \ra \Mor^\cA(M_0, M_2),
\end{align*}
we will view it as a map 
\begin{align*}
    \CFAh(-M_0)\boxtimes \CFDh(M_1) \otimes \Mor^\cA(M_1, M_2) \ra \CFAh(-M_0)\boxtimes \CFDh(M_1).
\end{align*}
Under this identification, fixing an element $f \in \Mor^\cA(M_1, M_2)$ yields a map 
\begin{align*}
    \CFAh(-M_0)\boxtimes \CFDh(M_1) \xra{\bI \boxtimes f} \CFAh(-M_0)\boxtimes \CFDh(M_1).
\end{align*}
This identification is unsurprising, but the correspondence between type A and D structures is nontrivial, so we justify this in \Cref{subsec:alg_lem}. 

\begin{proposition}\label{prop:FukayaHomSatellites_one_bspt}
	Let $C$ be a concordance from $K_0$ to $K_1$. Let $\btheta_0$ and $\btheta_1$ be the immersed curve invariants for $\KC_0$ and $\KC_1$ respectively. Identify $F_{S^3 \times I\smallsetminus C}$ as an element $f \in \hom_{\Fuk_\circ}(\btheta_0, \btheta_1)$.  
	Let $P \sub \bm{O}$ be a satellite pattern represented by a genus 1 bordered Heegaard diagram $\cH_P$ with beta curve $\bbeta_P$. Then, the following diagram
	\begin{align*}
		\begin{tikzcd}[ampersand replacement = \&]
			\hom_{\Fuk_{\circ\bullet}}(\btheta_0,\bbeta_P) \ar[r,"\simeq"] \ar[d,"m_2(f\text{,-})"] \& 
			\CFA^-(\cH_P) \boxtimes \CFDh(\KC_0)\ar[d,"\bI_P\boxtimes F_{S^3\times I \smallsetminus C}"]\\
			\hom_{\Fuk_{\circ\bullet}}(\btheta_1,\bbeta_P)  \ar[r,"\simeq"] \& 
			\CFA^-(\cH_P) \boxtimes \CFDh(\KC_1)
		\end{tikzcd}
	\end{align*}
	commutes up to homotopy. The horizontal arrows are given by Chen's pairing theorem.
\end{proposition}

\begin{rem}
	We expect that a generalization of \Cref{prop:FukayaHomSatellites_one_bspt} holds for the stronger invariant $\CFK_\cR$. There are a few issues in trying to apply the techniques of \cite{chen_hanselman_2023satellite}. Firstly, in Chen-Hanselman's work, they work with alpha-bordered Heegaard diagrams where one of the alpha circles is replaced with an immersed multicurve. Defining maps between such Heegaard diagrams necessitates working with immersed bordered triple diagrams with \emph{two} sets of alpha-arcs. In this context, the relevant moduli spaces appear to behave rather badly; these spaces do appear in \cite{LOT_bordered_HF}, but are not treated in full generality. Secondly, Chen-Hanselman make use of the fact that extensions of type D structures are unique up to homotopy equivalence. To apply their techniques, one would need an analogous statement for extensions of type D morphisms.
\end{rem}

\subsection{An algebraic lemma}\label{subsec:alg_lem}

It will be helpful to formulate Theorem 1 in terms of the tensor product version of the pairing theorem, rather than the $\Hom$ version. We therefore prove the following algebraic lemma.

Recall from \cite[Definition 2.3.16]{LOTBimodules2015} that, if $\cA$ is a dg-algebra, then its \emph{bar resolution} $\mathrm{Bar}(\cA)$ is the type-$\mathit{DD}$ bimodule whose underlying $\Bbbk$-module is the tensor algebra $T(\cA[1])$ and which has structure map $\delta^1$ given on basis elements by
\begin{align}
	\begin{split}
		\delta^1[a_1|\cdots|a_k]=&a_1\otimes[a_2|\cdots|a_k]\otimes 1+1\otimes[a_1|\cdots|a_{k-1}]\otimes a_k\\&+\sum_{1\leq j\leq k}1\otimes[a_1|\cdots|\mu_1(a_j)|\cdots|a_k]\otimes 1\\&+\sum_{1\leq j\leq k-1}1\otimes[a_1|\cdots|\mu_2(a_j,a_{j+1})|\cdots|a_k]\otimes 1.
	\end{split}
\end{align}
By \cite[Proposition 2.3.22]{LOTBimodules2015}, if $M$ is an operationally bounded right $A_\infty$-module over a dg-algebra $\cA$, then $M$ is $A_\infty$-homotopy equivalent to its bar resolution $\mathrm{Bar}(M):=M\boxtimes\mathrm{Bar}(\cA)\boxtimes\cA$ and $M\boxtimes\mathrm{Bar}(\cA)$ is isomorphic to $\overline{\overline{\mathrm{Bar}(\cA)}\boxtimes\overline{M}}$ by \cite[Lemma 2.15]{LOT_HF_as_morphism}. Since $\overline{\mathrm{Bar}(\cA)}\boxtimes\overline{M}$ is bounded whenever $M$ is, this tells us that any such $M$ is $A_\infty$-homotopy equivalent to a right dg-module of the form $\overline{N}_0\boxtimes\cA$ for some bounded left type-$D$ structure $N_0$.
\begin{definition}[{\cite[Definition 4.1]{CohenComposition}}]
	Given left type-$D$ structures $N_1$ and $N_2$ with preferred $\cA$-bases for their modulifications, we say a module homomorphism $f\in\Mor^\cA(N_1,N_2)$ is \emph{basic} if it is nonzero and there are generators $\bm{u}\in N_1$ and $\bm{v}\in N_2$, and an algebra element $\rho\in\cA$, such that $f:\bm{u}\mapsto\rho\bm{v}$ and $f$ vanishes on all other generators of $N_1$.
\end{definition}
\begin{proposition}\label{prop:PairingLemma}
	Let $(N_1,\delta_1^1)$ and $(N_2,\delta_2^1)$ be left type-$D$ structures over a dg-algebra $\cA$ over $\Bbbk$, and let $M$ be a bounded right $A_\infty$-module over $\cA$. Assume that $N_1$ and $N_2$ are homotopy equivalent to bounded type-$D$ structures. Then there is a homotopy commutative square
	\begin{align}
		\begin{tikzcd}[ampersand replacement=\&]
			M\boxtimes N_1\otimes\Mor^{\cA}(N_1,N_2)\arrow[r,"\id\boxtimes\mathit{ev}^1"]\arrow[d,"\simeq"'] \& M\boxtimes N_2\\
			\Mor^{\cA}(N_0,N_1)\otimes\Mor^{\cA}(N_1,N_2)\arrow[r,"\circ_2"] \& \Mor^{\cA}(N_0,N_2)\arrow[u,"\simeq"']
		\end{tikzcd},
	\end{align}
	where $\mathit{ev}^1:N_1\otimes\Mor^{\cA}(N_1,N_2)\to\cA\otimes N_2$ is the evaluation map, given by $\bm{x}\otimes f\mapsto f(\bm{x})$, and $\circ_2$ is the opposite composition map defined by
	\begin{align}
		\circ_2(f,g)=\begin{tikzcd}[ampersand replacement=\&,column sep=0.35cm]
			{} \& {}\arrow[d,dashed]\\
			{} \& f\arrow[d,dashed]\arrow[dl,bend right]\\
			\mu_2\arrow[d] \& g\arrow[d,dashed]\arrow[l]\\
			{} \& {}
		\end{tikzcd},
	\end{align}
	where we identify the complex of left module homomorphisms $\Mor^{\cA}(N_i,N_j)$ with the chain complex of $\Bbbk$-linear maps $f^1:N_i\to\cA\otimes N_j$ via the isomorphism sending each basic morphism $f:\bm{u}\mapsto\rho\bm{v}$ to the corresponding linear map $f^1:\bm{u}\mapsto\rho\otimes\bm{v}$, and $N_0$ is any left type-$D$ structure such that $M\simeq\overline{N}_0\boxtimes\cA$.
\end{proposition}
\begin{proof}
	The boundedness condition ensures that the box tensor products appearing in the square are well-defined. Note that if $M$ and $M'$ are homotopy equivalent $A_\infty$-modules, then there is a homotopy commutative square
	\begin{align}
		\begin{tikzcd}[ampersand replacement=\&]
			M\boxtimes N_1\otimes\Mor^{\cA}(N_1,N_2)\arrow[r,"\id\boxtimes\mathit{ev}^1"]\arrow[d,"\simeq"'] \& M\boxtimes N_2\\
			M'\boxtimes N_1\otimes\Mor^{\cA}(N_1,N_2)\arrow[r,"\id\boxtimes\mathit{ev}^1"] \& M\boxtimes N_2\arrow[u,"\simeq"']
		\end{tikzcd}
	\end{align}
	so we may assume without loss of generality that $M$ is of the form $M=\overline{N}_0\boxtimes\cA$. We then automatically have a commutative square
	\begin{align}
		\begin{tikzcd}[ampersand replacement=\&]
			M\boxtimes N_1\otimes\Mor^{\cA}(N_1,N_2)\arrow[r,"\id\boxtimes\mathit{ev}^1"]\arrow[d,equal] \& M\boxtimes N_2\\
			\overline{N}_0\boxtimes\cA\boxtimes N_1\otimes\Mor^{\cA}(N_1,N_2)\arrow[r,"\id\boxtimes\mathit{ev}^1"] \& \overline{N}_0\boxtimes\cA\boxtimes N_2\arrow[u,equal]
		\end{tikzcd}
	\end{align}
	and we claim that the square
	\begin{align}
		\begin{tikzcd}[ampersand replacement=\&]
			\overline{N}_0\boxtimes\cA\boxtimes N_1\otimes\Mor^{\cA}(N_1,N_2)\arrow[r,"\id\boxtimes\mathit{ev}^1"]\arrow[d,"\cong"'] \& \overline{N}_0\boxtimes\cA\boxtimes N_2\\
			\Mor^{\cA}(N_0,N_1)\otimes\Mor^{\cA}(N_1,N_2)\arrow[r,"\circ_2"] \& \Mor^{\cA}(N_0,N_2)\arrow[u,"\cong"']
		\end{tikzcd}
	\end{align}
	commutes, where the vertical isomorphisms $\overline{N}_i\boxtimes\cA\boxtimes N_j\cong\Mor^{\cA}(N_i,N_j)$ are the ones identifying a basic morphism $f:\bm{u}\mapsto\rho\bm{v}$ with the box tensor product element $\overline{\bm{u}}\boxtimes\rho\boxtimes\bm{v}$, where $\overline{\bm{u}}$ is the dual basis element of $\bm{u}$. In the graphical notation of \cite{LOT_bordered_HF}, the map $\id\boxtimes\mathit{ev}^1$ in the top row of this diagram is given by
	\begin{align}
		\begin{tikzcd}[ampersand replacement=\&,column sep=0.35cm]
			{} \& {}\arrow[d] \& {}\arrow[d,dashed] \& {}\arrow[ddl,bend left=20,rightsquigarrow]\\
			{} \& \otimes\arrow[d,Rightarrow] \& \delta_1\arrow[l,Rightarrow]\arrow[d,dashed] \& {}\\
			{} \& \otimes\arrow[d,Rightarrow] \& (-)\arrow[l]\arrow[d,dashed] \& {}\\
			\overline{\delta}_0\arrow[uuu,dashed]\arrow[r,Rightarrow] \& \mu\arrow[d] \& \delta_2\arrow[l,Rightarrow]\arrow[d,dashed] \& {}\\
			{}\arrow[u,dashed] \& {} \& {} \& {}
		\end{tikzcd},
	\end{align}
	where we use a squiggly arrow $\rightsquigarrow$ to denote elements of $\Mor^{\cA}(N_1,N_2)$, in the sense that
	\begin{align}
		\id\boxtimes\mathit{ev}^1(-\boxtimes(-\otimes f))=\begin{tikzcd}[ampersand replacement=\&,column sep=0.35cm]
			{} \& {}\arrow[d] \& {}\arrow[d,dashed] \& f\arrow[ddl,bend left=20,rightsquigarrow]\\
			{} \& \otimes\arrow[d,Rightarrow] \& \delta_1\arrow[l,Rightarrow]\arrow[d,dashed] \& {}\\
			{} \& \otimes\arrow[d,Rightarrow] \& (-)\arrow[l]\arrow[d,dashed] \& {}\\
			\overline{\delta}_0\arrow[uuu,dashed]\arrow[r,Rightarrow] \& \mu\arrow[d] \& \delta_2\arrow[l,Rightarrow]\arrow[d,dashed] \& {}\\
			{}\arrow[u,dashed] \& {} \& {} \& {}
		\end{tikzcd}=\begin{tikzcd}[ampersand replacement=\&,column sep=0.35cm]
			{} \& {}\arrow[d] \& {}\arrow[d,dashed] \& {}\\
			{} \& \otimes\arrow[d,Rightarrow] \& \delta_1\arrow[l,Rightarrow]\arrow[d,dashed] \& {}\\
			{} \& \otimes\arrow[d,Rightarrow] \& f\arrow[l]\arrow[d,dashed] \& {}\\
			\overline{\delta}_0\arrow[uuu,dashed]\arrow[r,Rightarrow] \& \mu\arrow[d] \& \delta_2\arrow[l,Rightarrow]\arrow[d,dashed] \& {}\\
			{}\arrow[u,dashed] \& {} \& {} \& {}
		\end{tikzcd}
	\end{align}
	as a linear map $M\boxtimes N_1\to M\boxtimes N_2$. As is typical in graphical calculi, upward arrows represent dual modules of arrows which are drawn downward. However, $\cA$ is a dg-algebra so we have
	\begin{align}
		\begin{tikzcd}[ampersand replacement=\&,column sep=0.35cm]
			{} \& {}\arrow[d] \& {}\arrow[d,dashed] \& {}\arrow[ddl,bend left=20,rightsquigarrow]\\
			{} \& \otimes\arrow[d,Rightarrow] \& \delta_1\arrow[l,Rightarrow]\arrow[d,dashed] \& {}\\
			{} \& \otimes\arrow[d,Rightarrow] \& (-)\arrow[l]\arrow[d,dashed] \& {}\\
			\overline{\delta}_0\arrow[uuu,dashed]\arrow[r,Rightarrow] \& \mu\arrow[d] \& \delta_2\arrow[l,Rightarrow]\arrow[d,dashed] \& {}\\
			{}\arrow[u,dashed] \& {} \& {} \& {}
		\end{tikzcd}=\begin{tikzcd}[ampersand replacement=\&,column sep=0.35cm]
			{} \& {}\arrow[d] \& {}\arrow[d,dashed] \& {}\arrow[dl,bend left=20,rightsquigarrow]\\
			{} \& \mu_2\arrow[d] \& (-)\arrow[l]\arrow[d,dashed] \& {}\\
			{}\arrow[uu,dashed] \& {} \& {} \& {}
		\end{tikzcd},
	\end{align}
	since $\mu=\mu_1+\mu_2$ in this case and there cannot be a $\mu_1$ term since $\mu$ automatically has at least two inputs --- coming from the $\cA$ and $\Mor^{\cA}(N_1,N_2)$ terms --- in $\id\boxtimes\mathit{ev}^1$. In other words, if $g$ is basic with $g:\bm{v}\to\sigma\bm{w}$ and $\overline{\bm{u}}\boxtimes\rho\boxtimes\bm{v}\in\overline{N}_0\boxtimes\cA\boxtimes N_1$ is nonzero, then
	\begin{align}
		(\id\boxtimes\mathit{ev}^1)(\overline{\bm{u}}\boxtimes\rho\boxtimes\bm{v}\otimes g)=\overline{\bm{u}}\boxtimes\mu_2(\rho,\sigma)\boxtimes\bm{w}=\overline{\bm{u}}\boxtimes\rho\sigma\boxtimes\bm{w},
	\end{align}
	which is identified with the basic morphism $\bm{u}\mapsto\rho\sigma\bm{w}$ under the isomorphism $\overline{N}_0\boxtimes\cA\boxtimes N_2\cong\Mor^{\cA}(N_0,N_2)$. This basic morphism is precisely the composition $\circ_2(f,g)$, where $f$ is the basic morphism $\bm{u}\mapsto\rho\bm{v}$ identified with $\overline{\bm{u}}\boxtimes\rho\boxtimes\bm{v}$ under the isomorphism $\overline{N}_0\boxtimes\cA\boxtimes N_1\cong\Mor^{\cA}(N_0,N_1)$. Since the basic morphisms form bases for $\Mor^{\cA}(N_i,N_j)$, this proves the desired result in the case that $N_0$ is bounded. However, if $N_0$ is not bounded, then it is homotopy equivalent to some bounded $N_0'$ and there is automatically a homotopy commutative square
	\begin{align}
		\begin{tikzcd}[ampersand replacement=\&]
			\Mor^{\cA}(N_0,N_1)\otimes\Mor^{\cA}(N_1,N_2)\arrow[r,"\circ_2"]\arrow[d,"\simeq"'] \& \Mor^{\cA}(N_0,N_2)\\
			\Mor^{\cA}(N_0',N_1)\otimes\Mor^{\cA}(N_1,N_2)\arrow[r,"\circ_2"] \& \Mor^{\cA}(N_0',N_2)\arrow[u,"\simeq"']
		\end{tikzcd}
	\end{align}
	so the bounded case implies the general case.
\end{proof}
\begin{corollary}
	If, for $i=1,2,3$, $Y_i$ is a bordered 3-manifold with boundary parametrized by a fixed pointed matched circle $\mathcal{Z}$, then there is a homotopy commutative square
	\begin{align}
		\begin{tikzcd}[ampersand replacement=\&]
			\CFAh(Y_1)\boxtimes\CFDh(Y_2)\otimes_\F\Mor^{\cA(-\mathcal{Z})}(Y_2,Y_3)\arrow[r,"\id\boxtimes\mathit{ev}^1"]\arrow[d,"\simeq"'] \& \CFAh(Y_1)\boxtimes\CFDh(Y_3)\\
			\CFh(Y_1\cup_\partial Y_2)\otimes_\F\CFh(-Y_2\cup_\partial Y_3)\arrow[r,"\widehat{F}_W"] \& \CFh(Y_1\cup_\partial Y_3)\arrow[u,"\simeq"']
		\end{tikzcd},
	\end{align}
	where $W:(Y_1\cup_\partial Y_2)\sqcup(-Y_2\cup_\partial Y_3)\to Y_1\cup_\partial Y_3$ is the pair-of-pants cobordism in the sense of \textup{\cite{CohenComposition}} and $\widehat{F}_W$ is the induced map on Heegaard Floer complexes.
\end{corollary}
\begin{proof}
	This is an immediate consequence of combining Proposition \ref{prop:PairingLemma} and \cite[Theorem 1.1]{CohenComposition}, which tells us that $\circ_2$ and $\widehat{F}_W$ fit into a homotopy commmutative square
	\begin{align}
		\begin{tikzcd}[ampersand replacement=\&]
			\Mor^{\cA(-\mathcal{Z})}(-Y_1,Y_2)\otimes_\F\Mor^{\cA(-\mathcal{Z})}(Y_2,Y_3)\arrow[r,"\circ_2"]\arrow[d,"\simeq"'] \& \Mor^{\cA(-\mathcal{Z})}(-Y_1,Y_3)\\
			\CFh(Y_1\cup_\partial Y_2)\otimes_\F\CFh(-Y_2\cup_\partial Y_3)\arrow[r,"\widehat{F}_W"] \& \CFh(Y_1\cup_\partial Y_3)\arrow[u,"\simeq"']
		\end{tikzcd}.
	\end{align}
\end{proof}

\begin{figure}
	\centering
	\includegraphics[width=.4\textwidth]{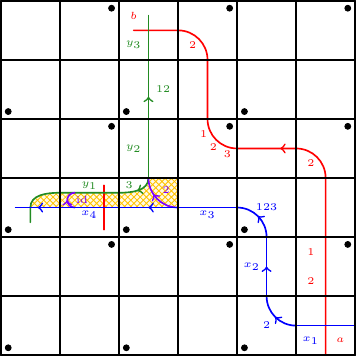}\hspace{1cm}\includegraphics[width=.4\textwidth]{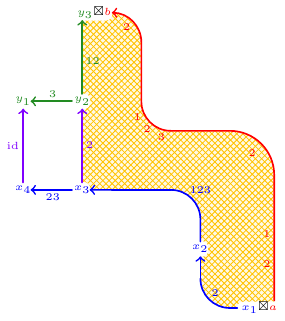}
	\caption{A triangle whose green boundary component cannot be consistently oriented by its constituent elementary arcs. We think of the intersection point between the blue and green curves as being a bounding chain.}
	\label{fig:TriangleCorrespondence}
\end{figure}

\begin{proof}[Proof of \Cref{prop:FukayaHomSatellites_one_bspt}]
    Embed $\btheta_{K_0}$ and $\btheta_{K_1}$ into the punctured torus according to the type D conventions. We will pair these curves with the beta curve of $\cH_P$. According to \cite{chen_hfk_satellite}, $\hom_{\Fuk_{\circ \bullet}}(\bbeta_P, \btheta_{K_i})$ computes $\CFK^-(P(K_i))$ for $i = 0,1$. By \Cref{cor:MorphismIdentification}, $\hom_{\circ}(\btheta_{K_0}, \btheta_{K_1})$ can be identified with $$\Mor^\cA(\CFDh(\KC_0), \CFDh(\KC_1))$$ by associating to an intersection point $\bm b_f$ of $\btheta_{K_0}\cap \btheta_{K_1}$ the morphism $f$ which counts triangles between $\btheta_{K_0}$, $\btheta_{K_1}$, and $\alpha_0\oplus\alpha_1$ with a corner at $\bm b_f$. We claim that counting triangles between $\btheta_{K_0}$, $\btheta_{K_1}$, and $\bbeta_P$ with a corner at $\bm b_f$ computes $\bI_P\boxtimes f.$

    Our argument is similar to the proof of \cite[Theorem 1]{chen_hfk_satellite}. In the pairing theorems of \cite{chen_hanselman_2023satellite} and \cite{hanselman2023bordered}, the type D structures are assumed to be (almost) reduced. Under this assumption, the two boundary components of any bigon contributing to the differential inherit orientations from the underlying immersed curves, since they are necessarily the unions of consistently oriented elementary arcs \Cref{fig:typeD_A_conventions}. These two oriented paths correspond to compatible type A and D operations, which pair up to produce terms in the box tensor product differential. For this argument, it is important that the curves are reduced. In our case, the analogue would be to assume that $f$ has no summands of the form $x \mapsto 1 \otimes y$. This is often not the case however, so generally there will indeed be triangles whose boundary components are made up of elementary arcs with inconsistent orientations. See \Cref{fig:TriangleCorrespondence}. 
	
	We can circumvent this issue by making use of the fact that we can explicitly pin down the morphism defined by $\bm b_f$. Let $D$ be a small ball around $\bm{b}_f$ and fix an identification of $(D, D\cap \btheta_0, D\cap \btheta_1)$ with $(D^2, [-1,1]\times \{0\}, \{0\} \times [-1,1])$ where $D^2$ is the unit disk in $\R^2$. Our convention is that any triangle emanating from $\bm{b}_f$ contains either the first or third quadrants of $D^2$. This partitions the set of triangles with a corner at $\bm{b}_f$ into two subsets, $\{T_i^1\}_{i=1}$ and $\{T_j^3\}_j$ which we call triangles of type 1 and type 3 respectively. Moreover, each subset is partially ordered by inclusion, so choose a labeling so that $T_{i}^k \subset T_{i+1}^k$. In particular, $T_0^1$ and $T_0^3$ are minimal. Each $T_i^k$ determines an element of $\cA$ according to its labeled corners on its $\alpha_0 \oplus \alpha_1$ boundary edge.
	
	Let $u:D^2\smallsetminus \{-i, 1, i\} \ra \Sigma$ be a holomorphic triangle and let $T$ be its image. To $T$ we associate a bigon $B_T$ as follows. Let $B_\varepsilon(1)$ be a ball of radius $\varepsilon$ centered at $1$ in $\C$. For small $\varepsilon$, $u(B_\varepsilon(1)\cap D^2)$ is contained in either $T_0^1$ or $T_0^3$, and we say that $T$ either type 1 or 3 accordingly. The two cases are symmetric, so assume $T$ is type 1. Choose the maximal $k$ such that $T_k^1 \cap \btheta_0 \sub T$ and $T_k^1 \cap \btheta_1 \sub T$. Let $\rho$ be the algebra element associated to $T_k^1$, and let $R_\rho$ be its associated elementary region in the sense of \cite[Figure 10]{chen_hfk_satellite}. Note that $T_k \smallsetminus T$ is precisely $R_\rho$. We can therefore define 
	\begin{align*}
		B(T) = (T\smallsetminus T_k)\cup R_\rho.
	\end{align*}
	Again, see \Cref{fig:TriangleCorrespondence} for an illustration. We will call the $R_\rho\cap \partial B_T$ part of the boundary its \emph{ghost edge}.
	
	With this terminology in place, we can identify the terms of $\bI_P\boxtimes f$ with triangles. First, suppose $T$ is a triangle from $x \otimes a$ to $y\otimes b$ furnishing a term in $m_1^{\bm{b}_f}$. Let $B_T$ be the associated bigon.  Just as in the proof of \cite[Theorem 1.2]{chen_hanselman_2023satellite}, applying the Chen collapse map determines an $A_\infty$ operation, 
	\begin{align*}
		m_{k+1}(a \otimes \rho_{L_1} \otimes \cdots \otimes \rho_{L_k}) = b,
	\end{align*}
	from $\bbeta_{P}\cap B_T$. The other half of the boundary consists of three pieces: $B_T \cap \btheta_0$, $B_T \cap \btheta_1$, and the ghost edge. The first determines a path $x = x_1 \xra{\rho_{I_1}} x_2 \xra{\rho_{I_2}} \cdots \xra{\rho_{I_{i-1}}} x_i$ in $\btheta_0$ and the second a path $y_1 \xra{\rho_{J_1}} y_2 \xra{\rho_{J_2}} \cdots \xra{\rho_{J_{j-1}}} y_j = y$ in $\btheta_1$. By construction of $B_T$, the ghost edge starts at $x_i$ and ends at $y_1$. In short, these determine a map
	\begin{align*}
		&(\bI_{\cA^{i+1}}\otimes \delta_1^{j})\circ (\bI_{\cA^{\otimes i}}\otimes (x_i \mapsto \rho \otimes y_1)) \circ \delta^i_0(x)\\&= \rho_{I_1}\otimes \cdots \rho_{I_i}\otimes \rho \otimes \rho_{J_1}\otimes \cdots \rho_{J_j}\otimes y.
	\end{align*}
	It follows from \cite{chen_hfk_satellite} that the Reeb chord sequences $\rho_{L_1} \otimes \cdots \otimes \rho_{L_k}$ and $\rho_{I_1}\otimes \cdots \rho_{I_i}\otimes \rho \otimes \rho_{J_1}\otimes \cdots \rho_{J_j}$ agree, giving rise to the term $x \otimes a\mapsto y \otimes b$ in $\bI_P \boxtimes f$.
	
	For the reverse direction, write $f = f_1 + \cdots + f_n$ as a sum of basic morphisms (i.e. those morphisms of the form $g(s) =\rho\otimes t.$) Equivalently, the map $f$ is determined by $n$ triangles in the torus; $f_\ell$ is the map induced by the $\ell$th triangle. Suppose $y \otimes b$ appears in $\bI_P\boxtimes f_\ell(x\otimes a)$, i.e. suppose that 
	\begin{align}\label{eqn:terms in Id box f}
		\begin{split}
			&\delta_0^i(x) = \rho_{I_1}\otimes \cdots \otimes \rho_{I_i} \otimes x_i \\
			&f_\ell(x_i) = \rho \otimes y_1 \\
			&\delta^j_1(y_1) = \rho_{J_1}\otimes \cdots \otimes\rho_{J_j}\otimes y\\
			&\mu_{2+i+j}(a\otimes \rho_{I_1}\otimes \cdots \otimes\rho_{I_i} \otimes\rho \otimes \rho_{J_1}\otimes \cdots\otimes \rho_{J_j}) = b.
		\end{split}
	\end{align}
	By adding a ghost edge from $x_i$ to $y_1$ labeled by $\rho$, this data determines a bigon, $B$, by the same argument as in \cite{chen_hfk_satellite}. Now, consider the portion of $\partial B$ containing the edge $x_i \mapsto \rho \otimes y_1$ produced by $f_\ell$ and let $R_\rho$ be its elementary region. Since terms of $f$ correspond to triangles emanating from $\bm{b}_f$, there is a triangle $T_\ell$ which corresponds to $f_\ell$, which can be used to define a triangle 
	\[
	T_B = (B\smallsetminus R_\rho)\cup T_\ell.
	\]
	This triangle clearly furnishes a term in $m_1^{\bm{b}_f}$ from $x\otimes a$ to $y \otimes b$ with corner at $\bm{b}_f = x_i \otimes y_1.$ This concludes the proof.
\end{proof}

	\section{Computations and Applications}\label{sec:applications}
In this section we consider a few applications of Theorems \ref{thm:FukayaHomPairing} and \ref{thm:FukayaHomSatellites}. First, we consider satellites of various disks, and the effectiveness of $\HFK$ in distinguishing them. We then consider a class of cobordisms, which we call \emph{splice cobordisms} built by gluing concordance complements together. We show that these cobordisms can naturally be computed using immersed curves and prove some detection results. Finally, we consider the bordered involutive invariants of \cite{HL_Inv_bordered_floer} and \cite{kang_bordered_involutive_HFK} and interpret those invariants in terms of immersed curves. 

\subsection{Satellite concordances via immersed curves}

Let $P$ be a pattern in the solid torus. The pattern $P$ induces an operator on the smooth (resp. topological) concordance group $\cC$ (resp. $\cC_{TOP}$), by taking the class $[K]$ to $[P(K)]$; this operation is well defined, since a concordance, $C$, from $K_0$ to $K_1$ has a satellite $P(C)$ from $P(K_0)$ to $P(K_1)$. In particular, these operators highlight the stark contrast in the smooth and topological categories. For instance, for $P$ an unknotted, winding number zero pattern, the operator 
\begin{align*}
	P: \cC_{TOP} \ra \cC_{TOP}
\end{align*}
is \emph{trivial}, since $P(K)$ is topologically slice, by the famous result of Freedman--Quinn \cite{freedman_quinn_top_4mflds}. On the other hand, in the smooth category, the kernel of such operators are completely mysterious. The most notable incarnation of this phenomenon concerns the \emph{Whitehead double pattern}, $\Wh$; conjecturally $\Wh$ is injective on the smooth concordance group (although it is trivial on the topological concordance group.) 

In \cite{guth2023doubled}, the authors considered an analogue of this phenomena for slice disks. When $P$ is an unknotted pattern (i.e., $P \ira \donut \ira S^3$ is the unknot), one can naturally define the notation of a satellite of a slice disk. Such a satellite pattern can be viewed as a map 
\begin{align*}
	P: \mathrm{Disk}_{\mathit{CAT}}(K) \ra \mathrm{Disk}_{\mathit{CAT}}(P(K))
\end{align*}
where $\mathrm{Disk}_{\mathit{CAT}}(K)$ denotes the set of smooth or topological isotopy (rel boundary) classes of slice disks for $K$, depending on whether $\mathit{CAT}$ is the smooth or topological category. Work of \cite{conway_powell_characterisation_homotopy_ribbon} implies that when $P$ has winding number zero, $P$ is trivial in the topological category; again, the kernels of such maps are completely mysterious in the smooth setting. \cite{guth2023doubled} gave some evidence that many unknotted satellites are injective in the smooth category by considering their effect on the maps on knot Floer homology.
\begin{defn}
	Let $P$ be an unknotted pattern in the solid torus. We say that $P$ is $\emph{$\HFKh$-injective}$ if for any slice disks $D_0$ and $D_1$ for $K$, $F_{D_0}\neq F_{D_1}$ implies that $F_{P(D_0)}\neq F_{P(D_1)}$.
\end{defn}

\begin{defn}\label{def:pos_slope}
	Let $\btheta_P$ be the immersed curve corresponding to an unknotted, $(1,1)$-pattern in the solid torus. We define the \emph{slope of $P$} to be the slope of the segment of $\btheta_P$ which crosses the curve $\btheta_U$ for the unknot after pulling $\btheta_P$ tight in the complement of the marked points in the universal cover.
\end{defn}

\begin{figure}
	\centering
	\includegraphics[scale=1]{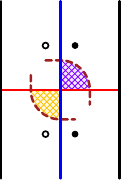}\hspace{1cm}\raisebox{1.5cm}{$\rightsquigarrow$}\hspace{1cm}\includegraphics[scale=1]{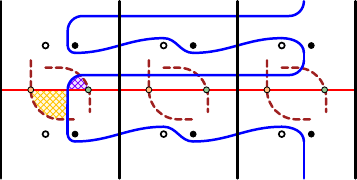}
	\caption{Left: Two possible triangles contributing terms to $f = F_{D_0} + F_{D_1}$. Right: The curve $\bbeta_\infty$ is replaced with $\bbeta_P$ which has positive slope.}
	\label{fig:cabled_disk_triangles}
\end{figure}

\begin{prop}
	Let $P$ be an unknotted $(1,1)$-pattern with positive slope. Then, $P$ is $\HFKh$-injective.
\end{prop}
\begin{proof}
	Let $D_0$ and $D_1$ be a pair of $\HFKh$-distinguishable slice disks for $K$. The map $f:=F_{D_0} + F_{D_1}$ can be realized as a linear combination of intersection points $\bm b_f$ between $\btheta_U$ and $\btheta_K$. Of course, these intersection points must lie on the $x$-axis. Since, $D_0$ and $D_1$ are distinguishable by $\HFKh$, there must be some triangle between $\bbeta_\infty$, $\btheta_U$, and $\btheta_K$. The map induced by a slide disk does not shift the Alexander grading; since the map on $\HFKh$ does count triangles which do not cover any of the basepoints, all three corners of this triangle must be contained in the region $S^1 \times [-\frac{1}{2},\frac{1}{2}]$; since the morphism $\bm b_f$ lies on the $x$-axis and the other two intersection points have Alexander grading zero. This greatly constrains the possible triangles contributing to $f$. In fact, there are only two possibilities, as shown in Figure \ref{fig:cabled_disk_triangles}; the boundary conditions require that the triangle emanate into the 1st or 3rd quadrants and must be completed by an arc of $\btheta_K$ which runs from the x-axis to the y-axis contained in the region $S^1 \times [-\frac{1}{2},\frac{1}{2}]$. 
    
    Now, to compute the maps induced by $P(D_0)$ and $P(D_1)$, pair $\Gamma_{f}$ with $\bbeta_P$. By assumption, the arc of $\bbeta_P$ which crosses the $x$-axis in the universal cover has positive slope. This arc must intersect $\btheta_K$ near one of the components of $f$; since this arc of $\btheta_P$ has positive slope, there is a triangle with the correct orientation of its edges to furnish a term in $F_{P(D_0)}+F_{P(D_1)}$. See Figure \ref{fig:cabled_disk_triangles}. Hence, the $(F_{P(D_0)}+F_{P(D_1)})(1) \in \HFKh(K)$ is nontrivial, and therefore the disks $P(D_0)$ and $P(D_1)$ cannot be smoothly isotopic rel boundary.  
\end{proof}

We can also easily detect stably exotic disks.

\begin{theorem}[\cite{guth_one_not_enough_exotic_surfaces}] Let $D$ and $D'$ be the two exotic positron disks of \cite{Hayden_CorksCoversComplex}, and let $J$ be their common boundary. Their $(n,1)$ cables, $D_n$ and $D_n'$, are exotic as well, and have stabilization distance at least $n$.
\end{theorem}

\begin{figure}
	\centering
	\includegraphics[scale=0.875]{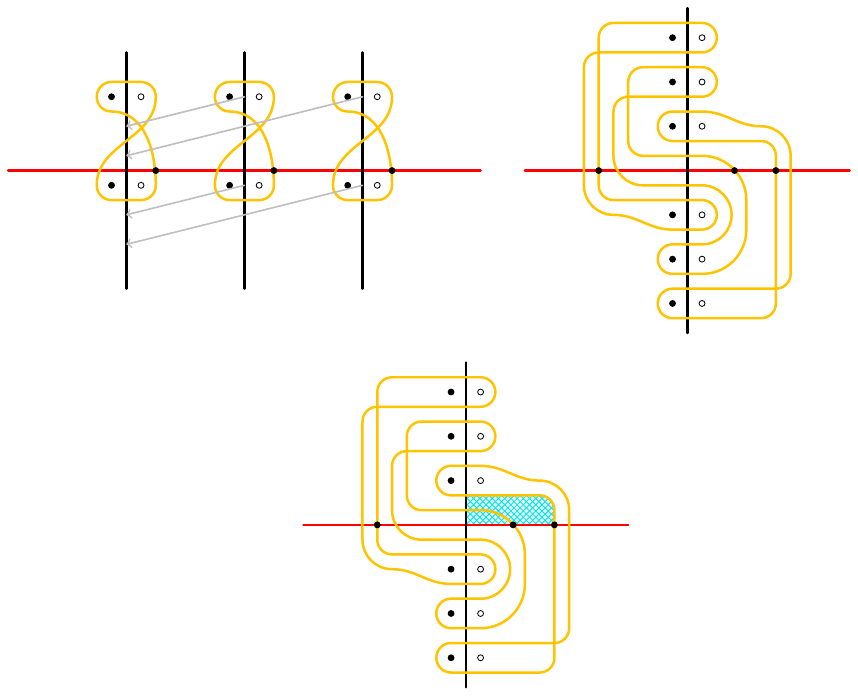}
	\caption{Applying the cabling procedure to the immersed curve for the pair $(D,D')$. }
	\label{fig:cabling_procedure}
\end{figure}

\begin{proof}
	The topological isotopy between $D$ and $D'$ produces an isotopy between the cabled disks, so it suffices to distinguish the disks smoothly. The immersed curve invariant for $J$, consists of four figure-8 components and a horizontal component which wraps around the cylinder. Two of the four figure 8 components do not contribute for grading reasons. In fact, the components of the immersed curve which will be relevant are exactly the same as the dashed red curves in \Cref{fig:9_46_disks}. 
    
    The disks $D$ and $D'$ induce maps 
	\begin{align*}
		F_D, F_{D'} \in \Hom_{\F[U]}(\CFK_{V=0}(U), \CFK_{V=0}(J))
	\end{align*}    
	By \Cref{thm:FukayaHomSatellites}, and taking $P$ to be the identity pattern, the morphism $f = F_D + F_{D'}$ corresponds to an element of $\hom_{\Fuk}(\btheta_U, \btheta_J)$, i.e. corresponds to an $\F$-linear combination of intersection points of $\btheta_U$ and $\btheta_J$, which we denote $\bm{b}_f$. The intersection points $\bm{b}_f$ agrees exactly with the green intersection points in \Cref{fig:9_46_disks}. By performing a change of basis, we can represent $f$ by a single intersection point between $\btheta_U$ and one of the figure 8 curves. There is then a single triangle contributing to the map $\bI \boxtimes f$.

	The map associated to the cabled disks can be obtained by applying the Hanselman-Watson cabling-for-immersed curve procedure for $(n, 1)$-cables to $(\btheta_U\cup \btheta_J, \bm{b}_f)$, or equivalently, by pairing with an immersed curve representative for the $(n, 1)$-cabling pattern, $\bbeta_n$; the $n=3$ case is shown in Figure \ref{fig:cabling_procedure}. There are now $n$ triangles which contribute to the sums of the cabled disk maps; consider the rightmost triangle in the diagram. This triangle terminates at a point $E$ in $\bbeta_{n}\cap \vartheta_J$. But note, there is a bigon which crosses the $V$-basepoint $n$ times and terminates at $E$ as well. Therefore, after setting $U = 0$, we see that $E$ descends to a generator of the homology with the property that $V^k \cdot E = 0$ only for $k \ge n$. Since multiplication by $V$ corresponds to a stabilization, it follows from \cite{juhasz_zemke_stabilization_bounds} that these disks can be stabilized $(n-1)$-times and still be distinguished by their induced maps on $\HFK^-(S^3,J)$.
\end{proof}

\subsection{Splicing concordances}
We now consider a class of examples which are quite well-suited to immersed curve computations. 

Let $X_0$ and $X_1$ be cobordisms with corners. We will divide the boundaries of such objects into horizontal and vertical parts: $\partial X_i = \partial_h X_i \cup \partial_v X_i$. Throughout, we will assume that $\partial_v X_i \cong T^2 \times I$, and equip $X_i$ with a parametrization $\varphi_i$ of $\partial_vX_i$. This data specifies a cobordism 
\begin{align*}
	X:= X_1 \cup_{\varphi_2 \circ \varphi_1^{-1}} X_2.
\end{align*}
Maps associated to cobordisms of this kind can be easily computed using bordered Floer homology.

\begin{proposition}
	Suppose that $X_1:Y_1^-\to Y_1^+$ and $X_2:Y_2^-\to Y_2^+$ are cobordisms with corners such that $\partial_vX_1\cong T^2 \times I\cong\partial_vX_2$ and let $X$ be as above. Then, there are maps $F_{X_1}\in \Mor_\cA(Y_1^-, Y_2^-)$ and $F_{X_2}\in \Mor^\cA(Y_1^+, Y_2^+)$ such that the diagram 
	\begin{align*}
		\begin{tikzcd}[ampersand replacement=\&]
			\CFAh(Y_1^-) \boxtimes \CFDh(Y_2^-) \ar[r,"\simeq"] \ar[d,"F_{X_1}\boxtimes F_{X_2}"] \& 
			\CFh(Y_1^- \cup Y_2^-) \ar[d,"F_X"]\\
			\CFAh(Y_1^+) \boxtimes \CFDh(Y_2^+) \ar[r,"\simeq"] \& 
			\CFh(Y_1^+ \cup Y_2^+)  
		\end{tikzcd}
	\end{align*}
	commutes up to homotopy.
\end{proposition}
\begin{proof}
	For $i=1,2$, fix a Morse function $f_i$ on $X_i$ such that $f_i|_{\partial X_i}$ has no critical points, as well as gradient like vector fields which are tangent to $\partial_vX_i$. These induce handle decompositions for $X_1$ and $X_2$; importantly, the attaching circles are contained in the interior of $Y_1^-$ and $Y_2^-$, and therefore, the same data determines a handle decomposition for $X$. 
	
	For simplicity, we only address the case that $X_1$ and $X_2$ consist entirely of two-handles attached along framed links $\cL_1$ and $\cL_2$; the case of 1- and 3-handles is simpler, and can be dealt with in the same way as in \cite{guth_one_not_enough_exotic_surfaces}. Fix (ordinary) Heegaard diagrams for $Y_i^- \smallsetminus \cL_i$, $\cH^i = (\Sigma_i, \alpha_1^i, \hdots, \alpha^i_{g-1}, \beta_{k_i+1}^i, \hdots, \beta^i_g)$, where $k_i$ is number of components of $\cL_i$. Promote these to bordered Heegaard triples as follows: 
	\begin{enumerate}
		\item For each $\beta^i_j$, let $\gamma_j^i$ be the curve obtained by a small Hamiltonian isotopy.
		\item For each component $K_\ell^i$ of $\cL_i$, let $\gamma_\ell^i$ be a longitude of $K_\ell^i$ determined by its framing.
		\item The parametrizations of $\partial Y_i^-$ determine preferred meridians and longitudes in $\Sigma_i$, denoted $\mu_i$ and $\lambda_i$, which are geometrically dual. Puncture $\Sigma_i$ at the intersection of these curves; the resulting arcs, together with our $\alpha$, $\beta$, and $\gamma$ curves makes $\Sigma\smallsetminus\nu(\pt)$ into a bordered Heegaard triple diagram, $\cH_B^i$.
	\end{enumerate}
	By construction, $\cH^1_B\cup \cH^2_B$ is a bordered Heegaard triple subordinate to a bouquet for $\cL_1 \cup \cL_2$; in particular, counting holomorphic triangles computes the map induced by $X$. By the pairing theorem for triangles, \cite{LOT_spectral_seq_II}, this map is homotopic to the box-tensor product of the triangle counting maps given by the triples $\cH_B^1$ and $\cH_B^2.$
\end{proof}

The simplest class of examples comes from concordances. If $C: K_0 \ra K_1$ is a concordance, then 
\begin{align*}
	\partial(S^3 \times I \smallsetminus C) \cong (S^3 \smallsetminus K_0) \cup T^2 \times I \cup (S^3, \smallsetminus K_1).
\end{align*}
If $C': J_0 \ra J_1$ is another concordance, we can glue them together by a map $\phi \times \bI: \partial (S^3 \smallsetminus J_0)
\times I \ra \partial (S^3 \smallsetminus J_1)\times I,$ and denote the resulting cobordism $W(C, C', \phi)$. If we insist that $\phi: \partial (S^3\smallsetminus K_0) \ra \partial (S^3\smallsetminus J_0)$ swaps meridians and longitudes, then $\partial_\pm(W(C, C', \phi))$ is the classical splice of $K_i$ and $J_i$ for $i = 0,1$. In this case, we will simply write $\cS(C, C')$ for $W(C, C',\phi)$, without explicitly recording the gluing information.

We will restrict to the even simpler case that one of the concordances is a product. We also note that, in a classical splice, since the two pieces are glued by a map which switches meridian and longitude, for any knot $J$, the splice of $S^3 \smallsetminus J$ and $S^3 \smallsetminus U$ is diffeomorphic to $S^3$. Analogously, any concordance $C: U \ra K$ can be spliced with the product $J \times I\subset S^3\times I$ to yield a cobordism which will be denoted $\cS_{C, J}$, from $S^3 \ra \cS(J, K)$. If $D$ is the slice disk determined by $C$, we will usually write $\cS_{D,J}$ for the cobordism obtained by capping off $\cS_{C, J}$ with a 4-ball. To demonstrate the utility of the immersed curve pairing, prove the following.

\begin{prop}\label{prop:splicing knots}
	Let $D$ and $D'$ be a pair of slice disks for $K$ which are distinguished by their induced map on $\HFKh$. Then, for any knot $J$ with $\tau(J) > 0$, the cobordisms $\cS_{D,J}$ and $\cS_{D',J}$ are distinguishable by $\HFh$.
\end{prop}
\begin{proof}
	Let $\btheta_K$, $\btheta_J$, and $\btheta_U$ be the immersed curve invariants for $K$, $J$, and $U$ respectively. Let $\bm b$ be a collection of intersection points between $\btheta_U$ and $\btheta_K$ representing the map $F_D + F_{D'}$.  Since $D$ and $D'$ can be distinguished by $\HFKh$, there is at least one triangle between $\btheta_U$, $\btheta_K$, and the infinity slope curve, which we denote $\bm \beta_\infty$. 
	
	Now, we compute the maps associated with $\cS_{D,J}$ and $\cS_{D',J}$. We do so by superimposing $\btheta_J$ on top of $(\btheta_U \cup \btheta_K, \bm b)$ so that their respective meridians and longitudes are interchanged. See Figure \ref{fig:trefoil_splice}.
	\begin{figure}
		\centering
		\includegraphics[scale=1]{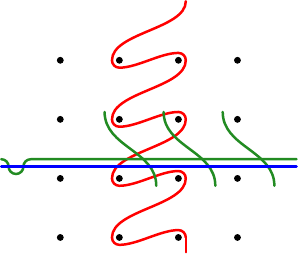}
		\caption{An example of \Cref{prop:splicing knots} in the case that $K$ is the trefoil.}
		\label{fig:trefoil_splice}
	\end{figure}    
	Note that this identification reflects $\btheta_J$. Since we assumed $\tau(J) > 0$, the distinguished component of $\btheta_J$ first intersects $\bm \beta_J$ at some positive height (in fact, at height $\tau(J)$); by the symmetry of the immersed curve invariants, it last intersects $\bm \beta_J$ at height $- \tau(J)$. Therefore, in the universal cover, after applying the splicing homeomorphism, there is an arc of positive slope between the points $(-\tau(J), -1)$ and $(\tau, 1).$ This line intersects $\btheta_U$ once, and this intersection point represents the generator of $\HFh(S^3).$ This line must also intersect $\btheta_K$ along an arc which formed one of the triangles distinguishing $D$ and $D'.$ It is then clear that there is a triangle furnishing a nonzero term in $F_{D_J} + F_{D_J'}$. See Figure \ref{fig:tau_splice} for an example.
	\begin{figure}
		\centering
		\includegraphics[scale=1]{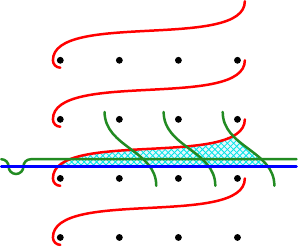}
		\caption{Triangles distinguishing the maps induced by the manifolds $\cS_{D,J}$ and $\cS_{D',J}$.}
		\label{fig:tau_splice}
	\end{figure}
\end{proof}

\subsection{$\SpinC$-Conjugation Action}

Finally, we consider bordered involutive Floer homology from the perspective of immersed curves. It is easy enough to translate results of \cite{HL_Inv_bordered_floer, kang_bordered_involutive_HFK} into the language of immersed curves. Though actually carrying out computations is quite tricky; as we shall see, explicitly identifying the domain and codomain is slightly subtle. In small examples, this is manageable, but this does seem to limit its utility. 

Let $M$ be a bordered 3-manifold with torus boundary. The $\SpinC$-conjugation action on the bordered modules are defined as maps
\begin{align}
	\iota^D_M:\CFDAh(\AZ) \boxtimes \CFDh(M) \ra \CFDh(M)
\end{align}
and
\begin{align}
	\iota^A_M: \CFAh(M)\boxtimes\CFDAh(\overline{\AZ})  \ra \CFAh(M).
\end{align}
Given two such bordered manifolds, $M$ and $N$ with bordered involutions $\iota^A_M$ and $\iota^D_N$, the involution $\iota_{M\cup N}$ can be computed as the box tensor product:
\begin{align}
	\iota^A_M\boxtimes \iota^D_N \simeq \iota_{M\cup N},
\end{align}
under the natural identification of $\CFAh(M)\boxtimes\CFDAh(\overline{\AZ}) \boxtimes \CFDAh(\AZ) \boxtimes \CFDh(N)$ with $\CFAh(M)\boxtimes \CFDh(N)$.
In particular, when $M$ is the exterior of a knot in $S^3$, the bordered involution $\iota_{\KC}$ is related to the action of $\iota_K$ on $\CFKh(S^3,K)$ by the following diagram: 
\begin{align}
	\begin{tikzcd}[ampersand replacement=\&]
		\CFAh(S^3 \smallsetminus K) \boxtimes \CFDh(T_\infty, \nu) \ar[r,"\simeq"] \ar[d,"\iota_{S^3 \smallsetminus K}^{-1}\boxtimes f"] \&
		\CFKh(S^3,K) \ar[dd,"\iota_K"] \\
		\overunderset{\hspace{-3mm}\displaystyle\CFAh(S^3 \smallsetminus K) \boxtimes \CFDAh(\overline{\AZ})}{\hspace{0.7mm}\displaystyle\CFDAh(\AZ)\boxtimes\CFDh(T_\infty, \nu)}{\boxtimes} \ar[d,"\simeq"] \& \\
		\CFAh(S^3 \smallsetminus K) \boxtimes \CFDh(T_\infty, \nu) \ar[r,"\simeq"]  \&
		\CFKh(S^3,K),
	\end{tikzcd}
\end{align}
which commutes up to homotopy for a particular choice of morphism
\begin{align}
	f \in \Mor^\cA(\CFDh(T_\infty, \nu), \CFDAh(\AZ) \boxtimes \CFDh(T_\infty, \nu)).
\end{align}
Unlike the 3-manifold case, the map $f$ is \emph{not} a homotopy equivalence; indeed, the two type D structures are not even homotopy equivalent. 

Understanding the 3-manifold story in terms of immersed curves is straightforward. Following \cite{hrw_properties_apps}, we can think of the functor $\CFDAh(\AZ) \boxtimes -$ as acting on immersed curves. We will denote by $\AZ(\btheta)$ the curves obtained from $\btheta$ by this action. It follows from \cite[Theorem 37]{hrw_properties_apps} that $\AZ$ acts trivially on the unlabeled immersed curves, and rotation by $\pi$ on labeled curves (the curves themselves are symmetric with respect to this symmetry). Therefore, by Theorem \ref{thm:FukayaHomPairing}, we may identify $\iota_M$ with a linear combination of intersection points of $\hom_{\Fuk_{\circ}}(\AZ(\btheta_M), \btheta_M),$ and compute the action of $\iota$ on glued manifolds by counting triangles. 

\begin{corollary}\label{cor:3mfoldIota}
	Let $M$ and $N$ be bordered manifolds with torus boundaries and bordered involutions $\iota_M$ and $\iota_N$. Let $\btheta_M$ and $\btheta_N$ be their immersed curve representatives.  Then, $\iota_{M}\boxtimes \iota_N$ can be computed as the composition
	\begin{align}
		\hom_{\Fuk_{\circ}}(\AZ(\btheta_M), \AZ(\btheta_N)) \xra{m_2(\iota_M,-)} \hom_{\Fuk_{\circ}}(\AZ(\btheta_M), \btheta_N) \xra{m_2(-,\iota_N)} \hom_{\Fuk_{\circ}}(\btheta_M, \btheta_N).
	\end{align}  
\end{corollary}

\begin{remark}
    There is a canonical identification of $\hom_{\Fuk_{\circ}}(\btheta_M, \btheta_N)$ with $\hom_{\Fuk_{\circ}}(\AZ(\btheta_M), \AZ(\btheta_N))$ due to the fact that there is a unique homotopy equivalence between  $\CFDAh(\overline{\AZ}\cup \AZ)$ and $\CFDAh(\bI)$ \cite[Corollary 8.4]{HL_Inv_bordered_floer}. Though, from the immersed curve perspective, it is unclear what this identification should be. In the examples we consider, we will try our best to avoid this issue; for our computations, it suffices to know that there is \emph{some} collection of points which determine the map. We will simply fix some identification, and use the fact that we know what the glued map is supposed to be in order to work backwards to determine which intersection points should contribute to the map. 
\end{remark}

Working with knots is more subtle. Presently, there is no Fukaya-categorical interpretation of morphisms between filtered type D structures, such as those arising from the modules for knots inside a bordered 3-manifold. In particular, it is unclear how to make sense of the map
\begin{align}
	f \in \Mor^{\cA(T^2)}(\CFDh(T_\infty, \nu), \CFDAh(\AZ) \boxtimes \CFDh(T_\infty, \nu))
\end{align} 
in terms of immersed curves. But, in the knot case, we really are working in a genus 1 bordered Heegaard diagram $\cH$ for the knot and the immersed curve invariant for $K$ is paired with this diagram. Therefore, our strategy will be simply realize $f$ as a map represented by counting triangles in a bordered Heegaard diagram.  

\begin{figure}
	\centering
	\includegraphics[width=.25\textwidth]{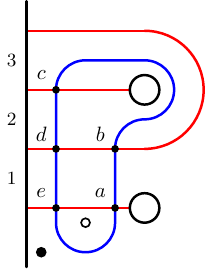}\hspace{1cm}\includegraphics[width=.325\textwidth]{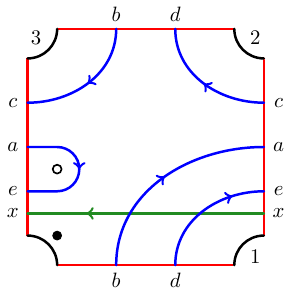}
	\caption{A bordered Heegaard diagram representing $\AZ\cup(T_\infty,\nu)$ (left), and the immersed curves $\bbeta_\infty$ and $\bbeta_{\AZ}$ in green and blue, respectively (right).}
	\label{fig:BorderedAZuTnu}
\end{figure}

\begin{lemma}[\cite{kang_bordered_involutive_HFK}]\label{lem:AZinfTorusTypeD}
	The bordered 3-manifold $\AZ \cup (T_\infty, \nu)$ is represented by the bordered Heegaard diagram in Figure \ref{fig:BorderedAZuTnu}. In particular, the type D structure $\CFDh(\AZ \cup (T_\infty, \nu))$ is given by 
	\begin{align}
		\begin{tikzcd}[ampersand replacement=\&]
			a \ar[d,dashed] \& 
			b \ar[l,"\rho_1" above] \ar[r,"\rho_3"] \ar[d,dashed] \&
			c \ar[dl, "\rho_2"] \\
			e \& 
			d. \ar[l,"\rho_1"] \& \\
		\end{tikzcd}
	\end{align}
	Further, there are two nontrivial homotopy classes of maps $\CFDh(T_\infty, \nu)\to\CFDh(\AZ\cup(T_\infty,\nu))$, realized by $f:=(x\mapsto a)$ and $g:=(x\mapsto e)$, respectively. The former is the map which recovers $\iota_K$.
\end{lemma}

Let $\bbeta_\infty$ and $\bbeta_{\AZ}$ be the immersed curves for $\CFDh(T_\infty, \nu)$ and $\CFDh(\AZ\cup(T_\infty, \nu))$, respectively (see Figure \ref{fig:BorderedAZuTnu}). The two intersection points between these curves \emph{do} realize the maps $f$ and $g$ (and are witnessed by the two small triangles in Figure \ref{fig:BorderedAZuTnu}). 

\begin{lemma}\label{lem:idF_triangles}
	Let $\btheta_K$ be the immersed curve invariant for a knot complement $\KC$. Then, the map $m_2(-,f)$ given by counting triangles with corner $f$ computes the map $\bI \boxtimes f$. 
\end{lemma}
\begin{figure}
	\centering
	\includegraphics[width = .7 \textwidth]{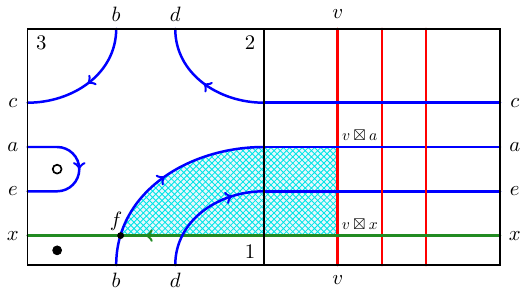}
	\caption{A triangle with corner at $f$ determined by $v\boxtimes x\in\theta\cap\bbeta_\infty$.}
	\label{fig:AZTriangle}
\end{figure}
\begin{proof}
	The model for $\CFAh(\KC)$ determined by $\btheta_K$ is reduced. Therefore, $\bI \boxtimes f$ acts on elements $v \in \iota_1\cdot \CFAh(\KC)$ by $v\boxtimes x \mapsto v \boxtimes a$,
	and trivially on $\iota_0\cdot \CFAh(\KC)$. Conversely, $\btheta_K$ is paired with $\bbeta_\infty$ and $\bbeta_{\AZ}$, it is clear that for every intersection point $v \boxtimes x$ of $\btheta_K \cap \bbeta_\infty$, there is a single triangle with corner at $f$; the third corner is $v \boxtimes a$. See Figure \ref{fig:AZTriangle}.
\end{proof}

\begin{remark}
	Curiously, every chain map from $\CFDh(\AZ\cup(T_\infty, \nu))$ to $\CFDh(T_\infty, \nu)$ is nullhomotopic. In particular, intersection points between $\bbeta_\infty$ and $\bbeta_{\AZ}$ cannot correspond to homotopy classes of chain maps.  This asymmetry is why it is necessary to tensor with $\iota_{\KC}^{-1}$ rather than $\iota_{\KC}$. 
\end{remark}

\begin{figure}
	\centering
	\includegraphics[width=.3\textwidth]{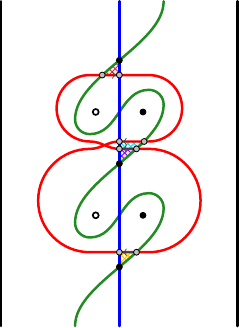}
	\caption{The nearest point map realizing $f$ in the lift of the curves $\theta$, $\bbeta_{\infty}$, and $\bbeta_{\AZ}$ to $S^1\times\R$. Note: for clarity, only one component of the lift of $\theta$ is shown.}
	\label{fig:Lift2}
\end{figure}

\begin{corollary}\label{cor:f_is_triangle}
	Let $\iota_\KC$ be the bordered involution for $\KC$. Then, the map 
    \begin{align*}
        \CFAh(S^3 \smallsetminus K) \boxtimes\CFDh(T_\infty, \nu) \xra{\iota_{\KC}^{-1}\boxtimes f} \CFAh(S^3 \smallsetminus K) \boxtimes \CFDAh(\overline{\AZ})\CFDAh(\AZ)\boxtimes\CFDh(T_\infty, \nu)
    \end{align*}
    as the composition: 
    \begin{align*}
        \hom_{\Fuk_{\circ\circ}}(\btheta_K, \bbeta_\infty) \xra{m_2(\iota_{\KC},-)} \hom_{\Fuk_{\circ\circ}}(\AZ(\btheta_K), \bbeta_\infty) \xra{m_2(-,f)} \hom_{\Fuk_{\circ\circ}}(\AZ(\btheta_K), \AZ(\bbeta_\infty)).
    \end{align*}
\end{corollary}
We note that after pulling tight in the infinite cylinder, $m_2(-,f)$ is just the nearest point map. Enjoy \Cref{fig:Lift2} and \Cref{fig:Lift}. 

\begin{example}
	Consider the $K = 4_1$. Its curve invariant $\btheta_K$ consists of the $x$-axis and a figure 8 component; it is drawn in orange in \Cref{fig:brieskorn_iota}. By pairing with $\bbeta_\infty$, we may read off the corresponding complex, which is shown in \Cref{fig:brieskorn_iota}. We have labeled the intersection points. 

    The curve invariant for $\AZ(\btheta_K)$ is identical to that of $\btheta_K$. Similarly for $\AZ(\bbeta_\infty)$. But, they are distinct as labeled curves, as $\AZ$ acts on curves by rotation by $180^\circ$. For this reason, when we pair $\AZ(\btheta_K)$ with $\AZ(\bbeta_\infty)$, the complex is almost identical, but we rotate all the labels. This fixes a (skew-equivariant) identification of $\hom_{\Fuk_{\circ\circ}}(\AZ(\btheta_K), \AZ(\bbeta_\infty))$ with $\hom_{\Fuk_{\circ\circ}}(\btheta_K, \bbeta_\infty).$ Let $\cF$ be the skew-equivariant morphism given by reflecting $\CFK_\cR(K)$ across the line $y = x$. Since $\btheta_K$ and $\AZ(\btheta_K)$ are identical as curves, we must add two additional intersections to make them admissible. The higher graded of these two intersection points is distinguished; when we pair with $\bbeta_\infty$, counting triangles with this intersection point at a corner is the nearest point map. In our case, according to our labeling, this intersection point realizes the flip map $\cF$.
    
    The action of $\iota_K$ was computed in \cite{hendricks_manolescu_Invol}, utilizing the fact that $\iota_k^2$ is homotopic to the Sarkar map. In the basis shown, $\iota_k = \cF + (a \mapsto x) + (x \mapsto e).$ It is easy to see that there are two triangles furnishing the maps $(a \mapsto x)$ and $(x \mapsto e)$, so we mark the corresponding intersection points of $\btheta_K$ and $\AZ(\btheta_K)$. These intersection points must correspond to $\iota_{\KC}^{-1}$, since these are the only triangles which can recover $\iota_K$ consistent with our identification of $\hom_{\Fuk_{\circ\circ}}(\AZ(\btheta_K), \AZ(\bbeta_\infty))$ with $\hom_{\Fuk_{\circ\circ}}(\btheta_K, \bbeta_\infty).$
	
	Having realized $\iota_{\KC}$ as a collection of intersection points, we can try to pair with other curves to compute the action of $\iota$ on other 3-manifolds. This is particularly simple when $\bbeta_n$ is the invariant for the $n$-framed solid torus. In this case, in which case there is a single endomorphism of the corresponding bordered invariant up to homotopy: the identity. Therefore, $\iota_{S^3_n(K)}$ can be computed by pairing  $\bbeta_n$ with $\btheta_K$ and $\AZ(\btheta_K)$ and counting triangles which contain $\iota_\KC$. There are three intersection points in $\hom_{\Fuk_{\circ}}(\btheta_K, \bbeta_n)$, which we label $A$, $B$, and $C$. There are also of course three intersection points in $\hom_{\Fuk_{\circ}}(\AZ(\btheta_K), \AZ(\bbeta_n))$, but a priori, it is unclear how to label them! If we permit ourselves to appear to Hanselman's work, and count bigons which cover the basepoints, we are forced to label points of $\hom_{\Fuk_{\circ}}(\AZ(\btheta_K), \AZ(\bbeta_n))$ according to the closest intersection point of $\hom_{\Fuk_{\circ}}(\btheta_K, \bbeta_n)$, as in Figure \Cref{fig:brieskorn_iota}. By counting triangles, we compute that $\iota = \id + (B \mapsto A.)$
	
	\begin{figure}
		\centering
		\includegraphics[width = .8 \textwidth]{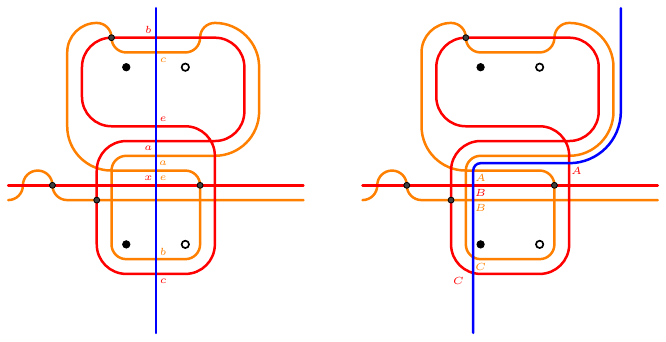}
		
		\begin{tikzcd}[column sep=0.5cm]
			a\ar[loop above,out=120,in=60,looseness=5]\ar[r] & x\ar[loop above,out=120,in=60,looseness=5]\ar[r] & e\ar[loop above,out=120,in=60,looseness=5]\\
			b\ar[r,bend right] & c\ar[l,bend right] &
		\end{tikzcd}\hspace{1.25cm}
		\begin{tikzcd}[column sep=0.5cm]
			A\ar[loop above,out=120,in=60,looseness=5] & B\ar[loop above,out=120,in=60,looseness=5]\ar[l] & C\ar[loop above,out=120,in=60,looseness=5]
		\end{tikzcd}
		\caption{The $\SpinC$-conjugation actions on $\HFKh(S^3, 4_1)$ and $\HFh(S^3_1(4_1))$ represented as morphisms between immersed curves.}
		\label{fig:brieskorn_iota}
	\end{figure}

\begin{example}
    As a final ``example'', we indulge ourselves by engaging in some speculation, and consider how one might formulate a ``immersed curve cabling formula for $\iota_K$''. By \cite{guth2024invariant}, the action of $\iota_{S^3\smallsetminus P(K)}$ is determined by $\iota_{\KC}$ as well as a \emph{type DA involution}
	\begin{align*}
		\iota_P: \cA \boxtimes \CFDAh(S^1\times D^2\smallsetminus P) \boxtimes \overline{\cA} \ra \CFDAh(S^1\times D^2\smallsetminus P).
	\end{align*}
	in the obvious way, $\iota_{S^3\smallsetminus P(K)} \sim \iota_P \boxtimes \iota_{\KC}$. When $P$ is a cabling pattern, the immersed curve invariant for $P(K)$ can be obtained by the Hanselman-Watson cabling formula (or also by the Chen-Hanselman formula). Roughly, the curve invariant $P(\btheta)$ for $P(K)$, is obtained from the invariant $\btheta$ for $K$ by choosing lifts of $\btheta$ to the universal cover, and then applying a certain diffeomorphism $f_P$ of the plane. The procedure defines an obvious map 
	\[
	\hom_{\Fuk_{\circ\circ}}(\AZ(\btheta), \btheta) \ra \hom_{\Fuk_{\circ\circ}}(P(\AZ(\btheta)), P(\btheta)),
	\]
	given by applying $\sum \bm{b} \mapsto f_P(\sum \bm{b})$. Counting triangles then yields \emph{some} endomorphism of $\CFKh(P(K))$. It seems natural to guess that this map is indeed $\iota_{P(K)}$. See Figure \ref{fig:figure8_cables_iota}. Note, that since this computation is purely speculative, we have also shaded disks that cross the basepoints, despite the fact that we have not proven anything about the action of $\iota_K$ on $\CFK_\cR(K)$.
	\begin{figure}
		\centering
		\includegraphics[scale=1]{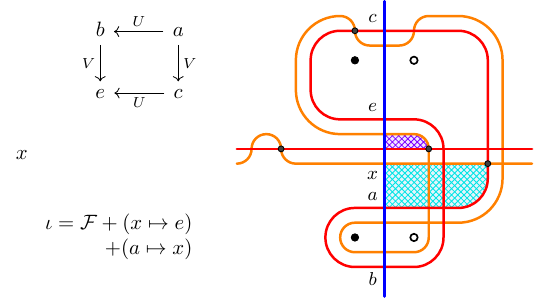}\\\vspace{0.25cm}\includegraphics[scale=1]{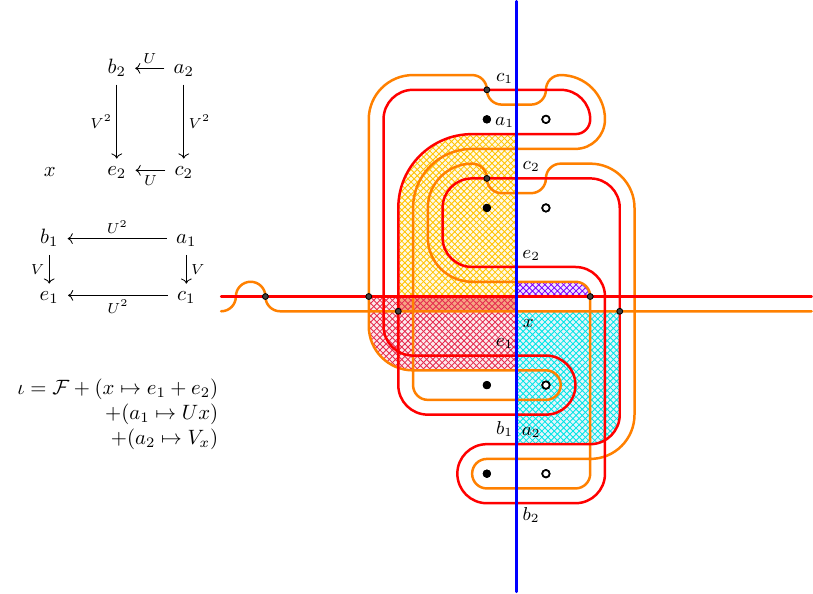}
		\caption{A endomorphism of $\CFK_\cR(K_{2,1})$ induced by $\iota_K$ and the diffeomorphism $f_P$.}
		\label{fig:figure8_cables_iota}
	\end{figure}
\end{example}
\end{example}

\begin{figure}
	\centering
	\includegraphics[width=0.8\textwidth]{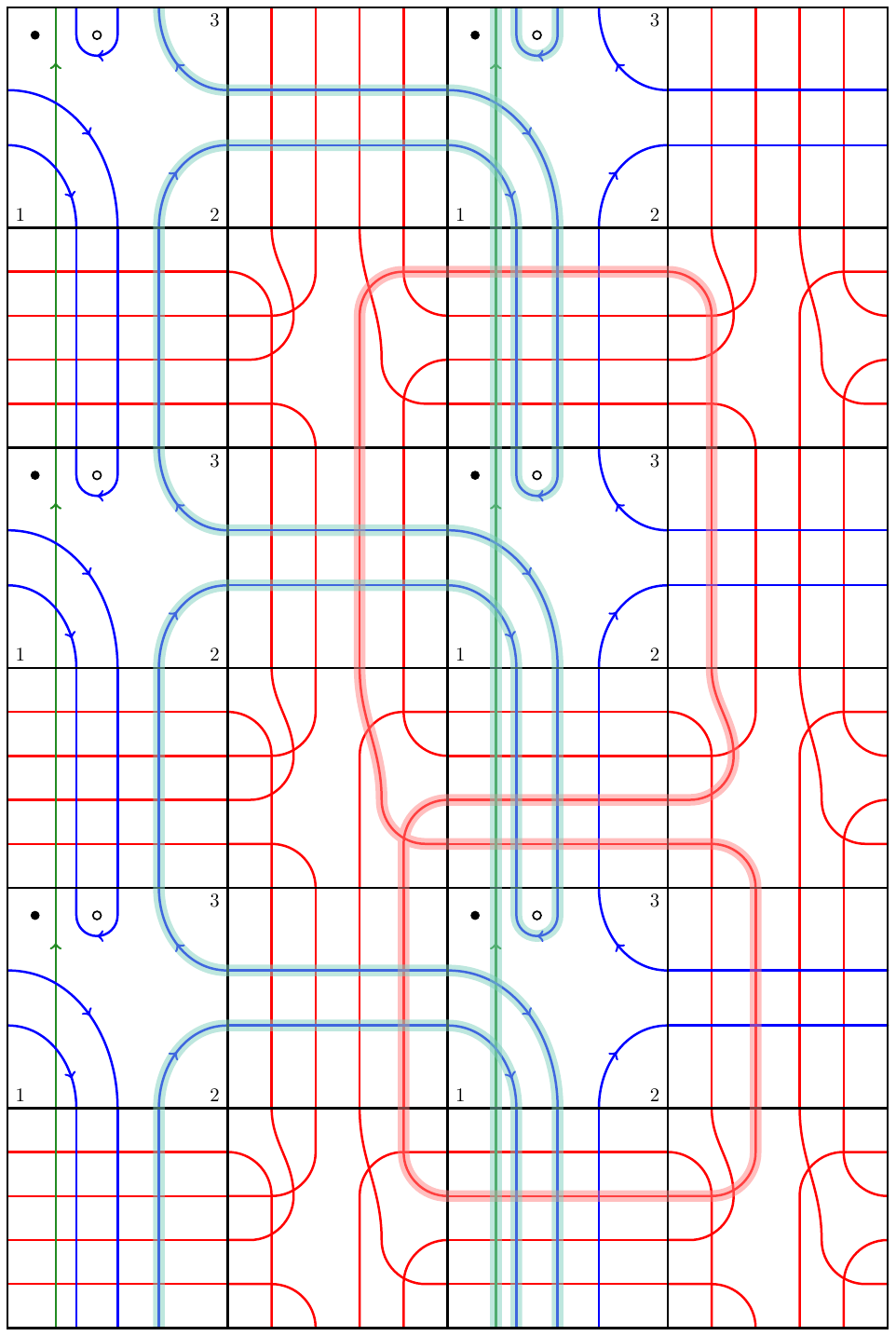}
	\caption{Lifts of the curves $\theta$ (red), $\bbeta_{\infty}$ (green), and $\bbeta_{\AZ}$ (blue) to the universal cover.}
	\label{fig:Lift}
\end{figure}

	\bibliographystyle{amsalpha}
	\bibliography{mathbib}
\end{document}